\documentclass[letterpaper,11pt]{amsart}
\title{Algebraic Geometry of  Gaussian Bayesian Networks}
\author{Seth Sullivant}
\address{Department of Mathematics and Society of Fellows, Harvard
  University, Cambridge, MA 02138}

\usepackage{latexsym,array,delarray,amsthm,amssymb,epsfig}

\theoremstyle{plain}
\newtheorem{thm}{Theorem}[section]
\newtheorem{lemma}[thm]{Lemma}
\newtheorem{prop}[thm]{Proposition}
\newtheorem{cor}[thm]{Corollary}
\newtheorem{conj}[thm]{Conjecture}

\theoremstyle{definition}
\newtheorem{defn}[thm]{Definition}
\newtheorem{ex}[thm]{Example}
\newtheorem{pr}[thm]{Problem}

\theoremstyle{remark}

\newcommand{\rr}{\mathbb{R}}
\newcommand{\cc}{\mathbb{C}}
\newcommand{\kk}{\mathbb{K}}
\newcommand{\bbe}{\mathbb{E}}

\newcommand{\bfa}{\mathbf{a}}

\newcommand{\bfx}{\mathbf{x}}
\newcommand{\bfy}{\mathbf{y}}
\newcommand{\bfz}{\mathbf{z}}

\newcommand{\cn}{\mathcal{N}}

\newcommand{\conv}{\mathrm{conv}}
\newcommand{\ind}{\mbox{$\perp \kern-5.5pt \perp$}}

\newcommand{\rmtop}{\mathrm{top}}
\newcommand{\rmpa}{\mathrm{pa}}

\begin{document}

\begin{abstract}
Conditional independence models in the Gaussian case are algebraic
varieties in the cone of positive definite covariance matrices.  We
study these varieties in the case of Bayesian networks, with a view
towards generalizing the recursive factorization theorem to situations
with hidden variables.  In the case when the underlying graph is a
tree, we show that the vanishing ideal of the model is generated by the conditional
independence statements implied by graph.  We also show that the
ideal of any Bayesian network is homogeneous with respect to a
multigrading induced by a collection of upstream random variables.
This has a number of important consequences for hidden variable
models.  Finally, we relate the ideals of Bayesian networks to a
number of classical constructions in algebraic geometry including
toric degenerations of the Grassmannian, matrix Schubert varieties,
and secant varieties. 
\end{abstract}

\maketitle


\section{Introduction}

A Bayesian network or directed graphical model is a statistical model
that uses a directed acyclic graph (DAG) to represent the conditional
independence structures between collections of random variables.  The
word Bayesian is used to describe these models because the nodes in
the graph can be used to represent random variables that
correspond to parameters or hyperparameters, though the basic
models themselves are not a priori Bayesian.  These models are used
throughout computational statistics to model complex interactions
between collections of random variables.  For instance, tree models
are used in computational biology for sequence alignment
\cite{Durbin1999} and in phylogenetics \cite{Felsenstein2004,
  Semple2003}.  Special cases of Bayesian networks include familiar
models from statistics like factor analysis \cite{Drton2006} and the
hidden Markov model \cite{Durbin1999}.

The DAG that specifies the Bayesian network specifies the model in two
ways.  The first is through a recursive factorization of the
parametrization, via restricted conditional distributions.  The second
method is via the conditional independence statements implied by the
graph.  The recursive factorization theorem \cite[Thm
  3.27]{Lauritzen1996} says that these two methods for specifying a
Bayesian network yield the same family of probability density
functions. 

When the underlying random variables are Gaussian or discrete,
conditional independence statements can be interpreted as algebraic
constraints on the parameter space of the global model.  In the
Gaussian case, this means that conditional independence corresponds to
algebraic constraints on the cone of positive definite matrices.  One
of our main goals in this paper is to explore the recursive
factorization theorem using algebraic techniques in the case of
Gaussian random variables, with a view  towards the
case of hidden random variables.  In this sense, the current paper is
a generalization of the work began in \cite{Drton2006} which concerned
the special case of factor analysis.  Some past work has been done
on the algebraic geometry of Bayesian networks in the discrete case in
\cite{Garcia2006, Garcia2005}, but there are many open questions that
remain in both the Gaussian and the discrete case. 

In the next section, we describe a combinatorial parametrization of a
Bayesian network in the Gaussian case.  In statistics, this
parametrization in known as the \emph{trek rule} \cite{Spirtes2000}.
We also describe the algebraic interpretation of
conditional independence in the Gaussian case which leads us to our
main problem:  comparing the vanishing ideal of the model $I_G$ to the
conditional independence ideal $C_G$.  Section \ref{sec:comp} describes the
results of computations regarding the ideals of Bayesian
networks, and some algebraic conjectures that these computations
suggest.  In particular, we conjecture that the coordinate ring of a
Bayesian network is always normal and
Cohen-Macaulay.   

As a first application of our algebraic perspective on Gaussian
Bayesian networks, we provide a new and greatly simplified proof of
the tetrad representation theorem  \cite[Thm 6.10]{Spirtes2000} in
Section \ref{sec:tetrad}.  Then in Section \ref{sec:tree} we provide
an extensive study of trees in the fully observed case.  In
particular, we prove that for any tree $T$, the ideal $I_T$ is a toric
ideal generated by linear forms and quadrics that correspond to
conditional independence statements implied by $T$.  Techniques from
polyhedral geometry are used to show that $ \cc[\Sigma]/I_T$ is always
normal and Cohen-Macaulay.

Sections \ref{sec:hidtree} and \ref{sec:alggeom} are concerned with
the study of hidden variable models.  In Section \ref{sec:hidtree} we
prove the Upstream Variables Theorem (Theorem \ref{thm:upstream}) which shows
that $I_G$ is homogeneous with respect to a two dimensional
multigrading induced by 
upstream random variables.  As a corollary, we deduce that hidden tree
models are generated by tetrad constraints.  Finally in Section
\ref{sec:alggeom} we show that models with hidden variables include,
as special cases, a number of classical constructions from algebraic
geometry.  These include toric degenerations of the Grassmannian,
matrix Schubert varieties, and secant varieties.

\subsection*{Acknowledgments}  I would like to thank Mathias Drton,
Thomas Richardson, Mike Stillman, and
Bernd Sturmfels for helpful comments and discussions about the results
in this paper.  The
IMA  provided funding and computer equipment while I worked on parts
of this project.


\section{Parametrization and Conditional Independence}

Let $G$ be a  directed acyclic graph (DAG) with vertex set $V(G) $ and
edge set $E(G)$.  Often, we will assume that $V(G) = [n] := \{1,2,
\ldots, n \}$.  To guarantee the acyclic assumption,    
we assume that the vertices are \emph{numerically ordered}; that is,
$i \to j  \in E(G)$ only if $i < j$. 
The Bayesian network associated to this graph can be specified by
either a recursive factorization formula or by conditional
independence statements.  We focus first on the recursive
factorization representation, and use it to derive an algebraic
description of the parametrization.  Then we introduce the
conditional independence constraints that vanish on the model and the
ideal that these constraints generate.

Let $X = (X_1, \ldots, X_n)$ be a random vector, and
let $f(x)$ denote the probability density function of this random
vector.  Bayes' theorem says that this joint density can be factorized
as a product 
$$f(x)  =  \prod_{i=1}^n  f_i(x_i | x_1, \ldots, x_{i-1}),$$
where $f_i(x_i | x_1, \ldots, x_{i-1})$ denotes the conditional density
of $X_i$ given $X_1 = x_1, \ldots, X_{i-1}= x_{i-1}$.  The recursive
factorization property of the graphical model is that each of the 
conditional densities $f_i(x_i | x_1, \ldots, x_{i-1})$ only depends on
the parents ${\rmpa}(i) = \{j \in [n] \,\, |\, \, j \to i \in E(G)
\}$.  We can rewrite this representation as 
$$f_i(x_i | x_1, \ldots, x_{i-1}) = f_i(x_i | x_{\rmpa(i)}).$$
Thus, a density function $f$ belongs to the Bayesian network if it
factorizes as
$$f(x)  =  \prod_{i =1}^n f_i(x_i | x_{\rmpa(i)}).$$

To explore the consequences of this parametrization in the Gaussian
case, we first need to recall some basic facts about Gaussian random
variables.  Each $n$-dimensional Gaussian random variable $X$ is
completely specified by its mean vector $\mu$ and its positive definite
covariance matrix $\Sigma$.  Given these data, the joint density
function of $X$ is given by
$$f(x) = \frac{1}{(2 \pi)^{n/2} |\Sigma|^{1/2}} \exp(- \frac{1}{2} (x -
\mu)^T \Sigma^{-1}(x - \mu) ),$$
where $|\Sigma|$ is the determinant of $\Sigma$.  Rather than writing
out the density every time, the shorthand $X \sim 
\mathcal{N}(\mu, \Sigma) $ is used to indicate that $X$ is a Gaussian
random variable with mean $\mu$ and covariance matrix $\Sigma$.  The
multivariate Gaussian 
generalizes the familiar ``bell curve'' of a univariate Gaussian and
is an important distribution in probability theory and multivariate
statistics because of the central limit theorem \cite{Bickel2001}.

Given an $n$-dimensional random variable $X$ and $A \subseteq [n]$,
let $X_A = (X_a)_{a \in A}$.  Similarly, if $x$ is a vector, then $x_A$ is
the subvector indexed by $A$.  For a matrix $\Sigma$, $\Sigma_{A,B}$
is the submatrix of $\Sigma$ with row index set $A$ and column index
set $B$.  Among the nice properties of Gaussian random variables are the fact
that marginalization and conditioning both preserve the Gaussian
property;  see \cite{Bickel2001}.

\begin{lemma}
Suppose that $X \sim \mathcal{N}(\mu, \Sigma)$ and let $A,B \subseteq
[n]$ be disjoint.  Then
\begin{enumerate}
\item  $X_A  \sim \mathcal{N}(\mu_A, \Sigma_{A,A})$ and
\item $X_A |  X_B = x_B  \sim  \mathcal{N}(\mu_A + \Sigma_{A,B}
  \Sigma_{B,B}^{-1}(x_B -\mu_B), \Sigma_{A,A} - \Sigma_{A,B}
  \Sigma_{B,B}^{-1} \Sigma_{B,A}).   $
\end{enumerate}

\end{lemma}

To build the Gaussian
Bayesian network associated to the DAG $G$, we allow any Gaussian
conditional distribution for the distribution $f(x_i | x_{\rmpa(i)})$.
This conditional distribution is recovered by saying that 
$$X_j  =  \sum_{i \in \rmpa(j)}  \lambda_{ij} X_i  +  W_j$$
where $W_j  \sim \cn(\nu_j, \psi_j^2)$ and is independent of the $X_i$
with $i < j$, and the $\lambda_{ij}$ are the regression parameters.
Linear transformations of Gaussian random variables are Gaussian, and thus $X$ is also a Gaussian random variable.
Since $X$ is 
completely specified by its mean $\mu$ and covariance matrix $\Sigma$,
we must calculate these from the conditional distribution.  The
recursive expression for the distribution of $X_j$ given the
variables preceding it yields a straightforward and recursive
expression for the mean and covariance.  Namely 
$$\mu_j  =  \bbe( X_j)  =  \bbe( \sum_{i \in \rmpa(j)}  \lambda_{ij}
X_i  +  W_j )  =  \sum_{i \in \rmpa(j)}  \lambda_{ij}  \mu_i  + \nu_j
$$ 
and if $k < j$ the covariance is:
\begin{eqnarray*}
\sigma_{kj} & =  &  \bbe\left((X_k - \mu_k)(X_j - \mu_j) \right) \\
 &  =  &  \bbe\left(   (X_k - \mu_k) \left(  \sum_{i \in \rmpa(j)}
 \lambda_{ij} (X_i - \mu_i) + W_j - \nu_j \right)  \right) \\ 
&  =   &  \sum_{i \in \rmpa(j)}  \lambda_{ij} \bbe\left( (X_k - \mu_k)
 (X_i - \mu_i) \right)  + \bbe \left( (X_k - \mu_k)(W_j - \nu_j)
 \right) \\ 
& =  &  \sum_{i \in \rmpa(j)}  \lambda_{ij}  \sigma_{ik}
\end{eqnarray*}
and the variance satisfies:
\begin{eqnarray*}
\sigma_{jj} & = & \bbe \left(( X_j - \mu_j)^2 \right) \\
&  = &  \bbe \left(  \left( \sum_{i \in \rmpa(j)}  \lambda_{ij} (X_i -
\mu_i) + W_j - \nu_j  \right)^2 \right) \\ 
& =  &  \sum_{i \in \rmpa(j)} \sum_{k \in \rmpa(j)}  \lambda_{ij}
\lambda_{kj} \sigma_{ik}  +  \psi^2_j. 
\end{eqnarray*}

If there are no constraints on the vector $\nu$,
there will be no constraints on $\mu$ either.  Thus, we will focus attention
on the constraints on the covariance matrix $\Sigma$.  If we further
assume that the $\psi^2_j$ are completely unconstrained, this will
imply that we can replace the messy expression for the covariance
$\sigma_{jj}$ by a simple new parameter $a_j$.  This leads us to the
algebraic representation of our model, called the \emph{trek rule} \cite{Spirtes2000}. 

For each edge $i\to j \in E(G)$ let $\lambda_{ij}$ be an indeterminate
and for each vertex $i \in V(G)$ let $a_i$ be an indeterminate.   Assume
that the vertices are numerically ordered, that is $i \to j \in E(G)$
only if $i < j$.  A \emph{collider} is a pair of edges $i \to k$,
$j \to k$ with the same head.   For each pair of vertices $i, j$, let
$T(i,j)$ be the collection of simple paths $P$ in $G$ from $i$ to $j$
such that there is no collider in $P$.  Such a colliderless path is
called a \emph{trek}.  The name trek come from the fact that every
colliderless path from $i$ to $j$ consists of a path from $i$ up to
some topmost element $\rmtop(P)$ and then from $\rmtop(P)$ back down
to $j$.  We think of each trek as a sequence of edges $k \to l$.   If $i = j$,
$T(i,i)$ consists of a single empty trek from $i$ to itself.

Let $\phi_G$ be the ring homomorphism
$$\phi_G:  \cc[\sigma_{ij} \, \, | \, \, 1 \leq i  \leq j \leq n]
\rightarrow  \cc[a_i, \lambda_{ij} \, \, | \, \, i,j \in [n] i \to j
  \in E(G) ]$$ 
$$\sigma_{ij}  \mapsto   \sum_{P \in T(i,j)}  a_{\rmtop(P)} \cdot
\prod_{ k \to l \in P }  \lambda_{kl}.$$ 
When $i = j$, we get $\sigma_{ii}  = a_i$.  If there is no trek in
$T(i,j)$, then $\phi_G(\sigma_{ij}) = 0$. 
Let $I_G = \ker \phi_G$.  Since $I_G$ is the kernel of a ring homomorphism, it is a prime ideal.

\begin{ex}\label{ex:4cycle}
Let $G$ be the directed graph on four vertices with edges $1 \to 2$,
$1 \to 3$, $2 \to 4$, and $3 \to 4$. 
The ring homomorphism $\phi_G$ is given by
$$
\begin{array}{cccc}
\sigma_{11}  \mapsto a_1 &   \sigma_{12} \mapsto  a_1 \lambda_{12} &
\sigma_{13} \mapsto a_1 \lambda_{13} & 
\sigma_{14} \mapsto a_1 \lambda_{12} \lambda_{24}  + a_1 \lambda_{13}
\lambda_{34}  \\ 
 &  \sigma_{22} \mapsto a_2  &  \sigma_{23} \mapsto a_1 \lambda_{12}
\lambda_{13}  &   
 \sigma_{24} \mapsto a_2 \lambda_{24} +  a_{1}  \lambda_{12}
 \lambda_{13} \lambda_{34}  \\ 
  &   &   \sigma_{33} \mapsto a_3  &   \sigma_{34} \mapsto a_3
 \lambda_{34} +  a_1 \lambda_{13} \lambda_{12}  
  \lambda_{24}  \\
  &   &   &  \sigma_{44} \mapsto a_4  \end{array}. $$
The ideal $I_G$ is the complete intersection of a quadric and a cubic:
$$I_G  =  \left<  \sigma_{11}\sigma_{23} - \sigma_{13} \sigma_{21},
\sigma_{12}\sigma_{23}\sigma_{34} + \sigma_{13} \sigma_{24}
\sigma_{23} + \sigma_{14}\sigma_{22} \sigma_{33} -
\sigma_{13}\sigma_{24} \sigma_{33} - \sigma_{13} \sigma_{22}
\sigma_{34} - \sigma_{14} \sigma_{23}^2  \right>. $$ 
\end{ex}

Dual to the ring homomorphism is the rational parametrization 
$$\phi^*_G :  \rr^{E(G) + V(G)}  \to  \rr^{ { n+1  \choose 2} }$$
$$  \phi^*_G(  a, \lambda)    =    (   \sum_{P \in T(i,j)}
a_{\rmtop(P)} \cdot \prod_{ k \to l \in P }  \lambda_{kl})_{i,j}.$$ 
We will often write $\sigma_{ij}(a, \lambda)$ to denote the coordinate
polynomial that represents this function. 

Let $\Omega  \subset \rr^{E(G) + V(G)}$ be the subset of parameter
space satisfying the constraints: 
$$a_i >  \sum_{ j \in \rmpa(i)}   \sum_{k \in  \rmpa(i)}  \lambda_{ji}
\lambda_{ki}  \sigma_{jk}(a, \lambda)$$ 
for all $i$, where in the case that $\rmpa(i)  = \emptyset$ the sum is zero.

\begin{prop}\label{prop:goodfac}[Trek Rule]
The set of covariance matrices in the Gaussian Bayesian network
associated to $G$ is the image $\phi^*_G(\Omega)$. 
In particular, $I_G$ is the vanishing ideal of the model.
\end{prop}

The proof of the trek rule parametrization
can also be found in \cite{Spirtes2000}. 

\begin{proof}
The proof goes by induction.  First, we make the substitution 
$$a_j =  \sum_{i \in \rmpa(j)} \sum_{k \in \rmpa(j)}  \lambda_{ij}
\lambda_{kj} \sigma_{ik}  +  \psi^2_j$$ 
which is valid because, given the $\lambda_{ij}$'s, $\psi^2_j$ can
be recovered from $a_j$ and vice versa.  Clearly $\sigma_{11} = a_1$.
By induction, suppose that the desired formula holds for all
$\sigma_{ij}$ with $i,j < n$.  We want to show that $\sigma_{in}$ has
the same formula.  Now from above, we have  
\begin{eqnarray*}
\sigma_{in}  & = & \sum_{k \in \rmpa(n)}  \lambda_{kn}  \sigma_{ik} \\
  &  =  &  \sum_{k \in \rmpa(n)}  \lambda_{kn}   \sum_{P \in T(i,k)}
  a_{\rmtop(P)} \cdot \prod_{ r \to s \in P }.  \lambda_{rs}   
\end{eqnarray*}
This last expression is a factorization of $\phi(\sigma_{kn})$ since every
trek from $i$ to $n$ is the union of a trek
$P\in T(i,k)$ and an edge $k \to n$ where $k$ is some parent of $n$. 
\end{proof}

The parameters used in the trek rule parametrization are a little
unusual because they involve a mix of the natural parameters
(regression coefficients $\lambda_{ij}$) and coordinates on the image
space (variance parameters $a_i$).  While this mix might seem unusual from a
statistical standpoint, we find that this
parametrization is rather useful for exploring the algebraic structure of the
covariance matrices that come from the model.  For instance: 

\begin{cor}
If $T$ is a tree, then $I_T$ is a toric ideal.
\end{cor}

\begin{proof}
For any pair of vertices $i,j$ in $T$, there is at most one trek
between $i$ and $j$.  Thus $\phi(\sigma_{ij})$
is a monomial and $I_T$ is a toric ideal. 
\end{proof}

In fact, as we will show in Section \ref{sec:tree}, when $T$ is a
tree, $I_T$ is generated by linear forms and quadratic binomials that
correspond to 
conditional independence statements implied by the graph.  Before
getting to properties of conditional independence, we first note that
these models are identifiable.  That is, it is possible to
recover the $\lambda_{ij}$ and $a_i$ parameters directly from $\Sigma$.
This also allows us to determine the most
basic invariant of $I_G$, namely its dimension.  

\begin{prop}\label{prop:dim}
The parametrization $\phi^*_G$ is birational.  In other
words, the model parameters $\lambda_{ij}$ and $a_i$ are identifiable
and $\dim I_G  =  \# V(G)  + \# E(G)$. 
\end{prop}

\begin{proof}
It suffices to prove that the parameters are identifiable via rational
functions of the entries of $\Sigma$, as all the
other statements follow from this.
We have $a_i = \sigma_{ii}$ so the $a_i$ parameters are identifiable.
We also know that for $i < j$
$$\sigma_{ij} = \sum_{k \in \rmpa(j)} \sigma_{ik} \lambda_{kj}.$$
Thus, we have the matrix equation
$$\Sigma_{\rmpa(j),j}  =  \Sigma_{\rmpa(j),\rmpa(j)}
\lambda_{\rmpa(j),j}$$
where $\lambda_{\rmpa(j),j}$ is the vector $(\lambda_{ij})^T_{i \in
  \rmpa(j)}$.  Since $\Sigma_{\rmpa(j), \rmpa(j)}$ is invertible in
the positive definite cone, we
have the rational formula
$$\lambda_{\rmpa(j),j} =
\Sigma_{\rmpa(j),\rmpa(j)}^{-1}\Sigma_{\rmpa(j),j}$$
and the $\lambda_{ij}$ parameters are identifiable.
\end{proof}

One of the
problems we want to explore is the connection between the prime ideal
defining the graphical model (and thus the image of the
parametrization) and the relationship to the ideal determined by the
independence statements induced by the model.  To explain this
connection, we need to recall some information about the algebraic
nature of conditional independence.  Recall the definition of
conditional independence.

\begin{defn}
Let $A$, $B$, and $C$ be disjoint subsets of $[n]$, indexing subsets
of the random vector $X$.  The conditional independence statement $A \ind B | C$ (``$A$ is independent of $B$ given $C$)
holds if and only if
$$f(x_A, x_B | x_C)  = f(x_A | x_C) f(x_B | x_C)$$
for all $x_C$ such that $f(x_C) \neq 0$.
\end{defn}

We refer to \cite{Lauritzen1996} for a more extensive introduction to
conditional independence.  In the Gaussian case, a conditional
independence statement is equivalent to an algebraic restriction on
the covariance matrix.

\begin{prop}
Let $A, B, C$ be disjoint subsets of $[n]$.  Then  $X \sim
\mathcal{N}(\mu, \Sigma)$ satisfies the conditional independence constraint
$A \ind B | C$ if and only if the submatrix 
$\Sigma_{A \cup C, B \cup C}$ has rank less than or equal to $\#C$.
\end{prop}

\begin{proof}
If $X  \sim  \mathcal{N}(\mu, \sigma)$, then 
$$X_{A \cup B}  |
X_C = x_C  \sim \mathcal{N}\left( \mu_{A\cup B} +  \Sigma_{A\cup B,C}
\Sigma_{C,C}^{-1}( x_C - \mu_C),   \Sigma_{A \cup B,A \cup B} -
\Sigma_{A\cup B,C} \Sigma_{C,C}^{-1}  \Sigma_{C,A \cup B} \right).$$  
The CI statement $A \ind B | C$ holds if and only if 
$(\Sigma_{A \cup B,A \cup B} - \Sigma_{A\cup B,C} \Sigma_{C,C}^{-1}
\Sigma_{C,A \cup B})_{A,B}  =  0$.  The $A,B$ submatrix of $\Sigma_{A
  \cup B,A \cup B} - \Sigma_{A\cup B,C} \Sigma_{C,C}^{-1}  \Sigma_{C,A
  \cup B}$ is easily seen to be $\Sigma_{A,B}  -  \Sigma_{A,C}
\Sigma_{C,C}^{-1}  \Sigma_{C,B}$ which is the Schur complement of the
matrix 
$$
\Sigma_{A \cup C,  B \cup C}  =  
\begin{pmatrix}
\Sigma_{A,B}  &   \Sigma_{A,C}  \\
\Sigma_{C,B}  &   \Sigma_{C,C}
\end{pmatrix}.
$$
Since $\Sigma_{C,C}$ is always invertible (it is positive definite),
the Schur complement is zero if and only if the matrix $\Sigma_{A \cup
  C,  B \cup C}$ has rank less than or equal to $\#C$. 
\end{proof}

Given a DAG $G$, a collection of conditional independence statements
are forced on the joint distribution by the nature of the graph.
These independence statements are usually described via the notion of
$d$-separation (the $d$ stands for ``directed'').   

\begin{defn}
Let $A$, $B$, and $C$ be disjoint subsets of $[n]$.  The set $C$
\emph{d-separates} $A$ and $B$ if every path in $G$ connecting a
vertex $i \in A$ and $B \in j$ contains a vertex $k$ that is either 
\begin{enumerate}
\item a non-collider that belongs to $C$ or
\item a collider that does not belong to $C$ and has no descendants
  that belong to $C$. 
\end{enumerate}
\end{defn}

Note that $C$ might be empty in the definition of $d$-separation.

\begin{prop}[\cite{Lauritzen1996}]
The conditional independence statement $A \ind B |C$ holds for the
Bayesian network associated to $G$ if and only if $C$ $d$-separates
$A$ from $B$ in $G$. 
\end{prop}

A joint probability distribution that satisfies all the conditional
independence statements implied by the graph $G$ is said to satisfy
the global Markov property of $G$.  The following theorem is a staple
of the literature of graphical models, that holds with respect to any
$\sigma$-algebra. 

\begin{thm}[Recursive Factorization Theorem] \cite[Thm 3.27]{Lauritzen1996}
A probability density has the recursive factorization property with
respect to $G$ if and only if it satisfies the global Markov
property. 
\end{thm}

\begin{defn}
Let $C_G \subseteq \cc[\Sigma]$ be the ideal generated by the minors
of $\Sigma$ corresponding  to the conditional independence statements
implied by $G$;  that is, 
$$C_G   =  \left<   (\#C+1) \mbox{ minors of }
\Sigma_{A \cup C,  B \cup C}  \,\, | \, \,   C  \mbox{ $d$-separates }
A \mbox{ from } B \mbox{ in } G  \right>.$$ 
The ideal $C_G$ is called the \emph{conditional independence ideal} of $G$.
\end{defn}

A direct geometric consequence of the recursive factorization theorem
is the following 

\begin{cor}
For any DAG $G$, 
$$V(I_G)   \cap  PD_n   =  V(C_G)  \cap PD_n.$$
\end{cor}

In the corollary  $PD_n  \subset  \rr^{ {n+1  \choose 2}}$ is the cone
of $n \times n$ positive definite symmetric matrices.  It seems
natural to ask whether or not $I_G  = C_G$ for all DAGs $G$.  For
instance, this was true for the DAG in Example \ref{ex:4cycle}.  The
Verma graph provides a natural counterexample. 

\begin{ex}\label{ex:verma}
Let $G$ be the DAG on five vertices with edges $1 \to 3$, $1 \to 5$,
$2 \to 3$, $2 \to 4$, $3 \to 4$, and $4 \to 5$.  This graph is often
called the Verma graph.   

\begin{figure}[h]
\begin{center}  \resizebox{5cm}{!}{\includegraphics{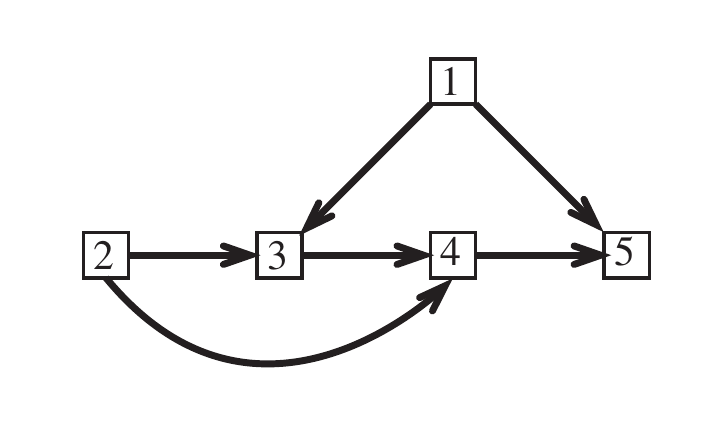}}
\end{center}
\end{figure}
The conditional independence statements implied by the model are all
implied by the three statements $1 \ind 2$, $1 \ind 4 | \{2,3\}$, and
$\{2,3\} \ind 5 | \{1,4\}$.  Thus, the conditional independence ideal
$C_G$ is generated by one linear form  and five
determinantal cubics.  In this case, we find that $I_G =  C_G + \left<
f \right>$ where $f$ is the degree four polynomial: 
\begin{eqnarray*}
f & = &
 \sigma_{23}\sigma_{24}\sigma_{25}\sigma_{34} -
 \sigma_{22}\sigma_{25}\sigma_{34}^2 -
 \sigma_{23}\sigma_{24}^2\sigma_{35} +
 \sigma_{22}\sigma_{24}\sigma_{34}\sigma_{35}\\   
 &   &    -\sigma_{23}^2\sigma_{25}\sigma_{44} +
 \sigma_{22}\sigma_{25}\sigma_{33}\sigma_{44}  +
 \sigma_{23}^2\sigma_{24}\sigma_{45} -
 \sigma_{22}\sigma_{24}\sigma_{33}\sigma_{45}. 
 \end{eqnarray*}
We found that the primary decomposition of $C_G$ is 
$$C_G \quad =  \quad I_G \cap \left< \sigma_{11}, \sigma_{12},
\sigma_{13}, \sigma_{14} \right>$$ 
so that $f$ is not even in the radical of $C_G$.  Thus, the zero set
of $C_G$ inside the positive \emph{semi}definite cone contains 
singular covariance matrices that are not limits of distributions that
belong to the model.  Note that since none of the indices of the
$\sigma_{ij}$ appearing in $f$ contain $1$, $f$ vanishes on the
marginal distribution for the random vector $(X_2,X_3,X_4,X_5)$.  This
is the Gaussian version of what is often called the Verma constraint.
Note that this computation shows that the Verma constraint is still
needed as a generator of the unmarginalized Verma model.  \qed 
\end{ex}

The rest of this paper is concerned with studying the ideals $I_G$ and
investigating the circumstances that guarantee that $C_G = I_G$.  We
report on results of a computational study in the next section.
Towards the end of the paper, we study the ideals $I_{G,O}$ that arise
when some of the random variables are hidden.


\section{Computational Study}\label{sec:comp}

Whenever approaching a new family of ideals, our first instinct is to
compute as many examples as possible to gain some intuition about the
structure of the ideals.  This section summarizes the results of our
computational explorations.

We used Macaulay2 \cite{M2} to compute the generating sets of all ideals $I_G$
for all DAGs $G$ on $n \leq 6$ vertices.  Our computational results
concerning the problem of when $C_G = I_G$ are summarized in the
following proposition.

\begin{prop}
All DAGs on $n \leq 4$ vertices satisfy $C_G = I_G$.  Of the $302$
DAGs on $n =5$ vertices, exactly $293$ satisfy $C_G = I_G$.  Of the
$5984$ DAGs on $n=6$ vertices exactly $4993$ satisfy $C_G = I_G$.
\end{prop}

On $n=5$ vertices, there were precisely nine graphs that fail to
satisfy $C_G = I_G$.  These nine exceptional graphs are listed below.
The numberings of the DAGs come from the Atlas of Graphs \cite{Read1998}.
Note that the Verma graph from Example \ref{ex:verma} appears as
$A_{218}$ after relabeling vertices.

\begin{enumerate}
\item $A_{139}$:  $1 \to 4$, $1 \to 5$, $2 \to 4$, $3 \to 4$, $4 \to
  5$.
\item $A_{146}$:  $1 \to 3$, $2 \to 3$, $2 \to 5$, $3 \to 4$, $4 \to 5$.
\item $A_{197}$:  $1 \to 2$, $1 \to 3$, $1 \to 5$, $2 \to 4$, $3 \to 4$, $4 \to 5$.
\item $A_{216}$: $1\to 2$, $1 \to 4$, $2 \to 3$, $2 \to 5$, $3 \to 4$, $4 \to 5$.
\item $A_{217}$:  $1 \to 3$, $1 \to 4$, $2 \to 4$, $2 \to 5$, $3 \to
  4$, $4 \to 5$.
\item $A_{218}$:  $1 \to 3$, $1 \to 4$, $2 \to 3$, $2 \to 5$, $3 \to 4$, $4 \to 5$.
\item $A_{275}$:  $1 \to 2$, $1 \to 4$, $1 \to 5$, $2 \to 3$, $2 \to 5$, $3 \to 4$, $4 \to 5$.
\item $A_{277}$:  $1 \to 2$, $1 \to 3$, $1 \to 5$, $2 \to 4$, $3 \to 4$, $3 \to 5$, $4 \to 5$.
\item $A_{292}$:  $1 \to 2$, $ 1\to 4$, $2 \to 3$, $2 \to 5$,   $3 \to 4$, $3 \to 5$, $4 \to 5$.
\end{enumerate} 

The table below displays the numbers of minimal generators of different degrees for each of the ideals $I_G$ where $G$ is one of the nine graphs on five vertices such that $C_G \neq I_G$.  The coincidences among rows in this table arise because sometimes two different graphs yield the same family of probability distributions.  This phenomenon is known as Markov equivalence \cite{Lauritzen1996, Spirtes2000}.

\smallskip

\begin{center}
\begin{tabular}{|c||c|c|c|c|c|}
\hline
Network & 1 & 2 & 3 & 4 & 5  \\  \hline \hline
$A_{139}$ &  3 & 1 & 2 & 0 & 0 \\ \hline
$A_{146}$ &  1 & 3 & 7 & 0 & 0 \\ \hline
$A_{197}$ &  0  & 1 & 5 & 0 & 1 \\  \hline
$A_{216}$ &  0 & 1 & 5 & 0 & 1 \\ \hline
$A_{217}$ & 2 & 1 & 2 & 0 & 0 \\ \hline
$A_{218}$ & 1 & 0 & 5 & 1 & 0 \\ \hline
$A_{275}$ & 0 & 1 & 1 & 1 & 3 \\ \hline
$A_{277}$ & 0 & 1 & 1 & 1 & 3 \\ \hline
$A_{292}$ & 0 & 1 & 1 & 1 & 3 \\ \hline
\end{tabular}
\end{center}

\smallskip

It is worth noting the methods that we used to perform our
computations, in particular, how we computed generators for the ideals
$I_G$.  Rather than using the trek rule directly, and computing the
vanishing ideal of the parametrization, we exploited the recursive
nature of the parametrization to determine $I_G$.  This is summarized
by the following proposition.

\begin{prop}
Let $G$ be a DAG and $G \setminus n$ the DAG with vertex $n$ removed.
Then
$$I_G  =  \left(  I_{G \setminus n}  +  \left<   \sigma_{in}  -
\sum_{j \in \rmpa(n)}  \lambda_{jn} \sigma_{ij}  \, \, | \, \, i \in
    [n-1] \right>  \right)   \bigcap \cc[\sigma_{ij} \, \, | \, \, i,j
      \in [n] ]$$
where the ideal $I_{G \setminus n}$ is considered as a graph on $n-1$
vertices. 
\end{prop}

\begin{proof}
This is a direct consequence of the trek rule:  every trek that goes
to $n$ passes through a parent of $n$ and cannot go below $n$.
\end{proof}

Based on our (limited) computations up to $n = 6$ we propose
some optimistic conjectures about the structures of the ideals $I_G$.

\begin{conj}\label{conj:sat}
$$I_G =  C_G :  \prod_{A \subset [n]}   (| \Sigma_{A,A}|)^{\infty}$$
\end{conj}

Conjecture \ref{conj:sat} says that all the uninteresting components
of $C_G$ (that is, the components that do not correspond to
probability density functions) lie on the boundary of the positive
definite cone.  Conjecture \ref{conj:sat} was verified for all DAGs on $n \leq 5$ vertices.  Our computational evidence also suggests that all the
ideals $I_G$ are Cohen-Macaulay and normal, even for graphs with loops
and other complicated graphical structures.  

\begin{conj}\label{conj:cm}
The quotient ring $\cc[\Sigma]/I_G$ is normal and Cohen-Macaulay for all $G$.
\end{conj}

Conjecture \ref{conj:cm} was verified computationally for all graphs
on $n \leq 5$ vertices and graphs with $n= 6$ vertices and less than
$8$ edges.
We prove Conjecture \ref{conj:cm} when the underlying graph is a tree
in Section \ref{sec:tree}.  A more negative conjecture concerns the
graphs such that $I_G = C_G$.

\begin{conj}
The proportion of DAGs on $n$ vertices such that $I_G = C_G$ tends to
zero as $n \to \infty$.
\end{conj}

To close the section, we provide a few useful propositions for
reducing the computation of the generating set of the ideal $I_G$ to
the ideals for smaller graphs.

\begin{prop}\label{prop:disconnect}
Suppose that $G$ is a disjoint union of two subgraph $G = G_1 \cup G_2$.  Then
$$I_G =  I_{G_1} + I_{G_2} + \left< \sigma_{ij} \, \, | \, \, i \in
V(G_1), j \in V(G_2) \right>.$$ 
\end{prop}

\begin{proof}
In the parametrization $\phi_G$, we have $\phi_G(\sigma_{ij}) = 0$ if
$i \in V(G_1)$ and $j \in V(G_2)$, because there is no trek from $i$ to $j$.
Furthermore,  $\phi_G(\sigma_{ij}) = \phi_{G_1}(\sigma_{ij})$ if
$i,j \in V(G_1)$ and $\phi_G(\sigma_{kl})  = \phi_{G_2}(\sigma_{kl})$
if $k,l \in V(G_2)$ and these polynomials are in disjoint sets of
variables.  Thus, there can be no nontrivial relations involving both
$\sigma_{ij}$ and $\sigma_{kl}$. 
\end{proof}

\begin{prop}
Let $G$ be a DAG with a vertex $m$ with no children and a
decomposition into two induced subgraphs $G = G_1 \cup G_2$ such that
$V(G_1) \cap V(G_2) = \{m\}$.  Then 
$$I_G = I_{G_1} + I_{G_2} + \left< \sigma_{ij} \, \, | \, \, i \in
V(G_1) \setminus \{m\},  j \in V(G_2) \setminus \{m\} \right>.$$
\end{prop}

\begin{proof}
In the paremtrization $\phi_G$, we have $\phi_G(\sigma_{ij}) = 0$ if
$i \in V(G_1) \setminus \{m\}$ and $j \in V(G_2) \setminus \{m\}$,
because there is no trek from $i$ to $j$. Furthermore
$\phi_G(\sigma_{ij}) = \phi_{G_1}(\sigma_{ij})$ if $i,j \in V(G_1)$
and $\phi_G(\sigma_{kl}) = \phi_{G_2}(\sigma_{kl})$ if $k,l \in
V(G_2)$ and these polynomials are in disjoint sets of variables unless
$i = j = k = l = m$.  However, in this final case,
$\phi_G(\sigma_{mm}) = a_m$ and this is the only occurrence of $a_m$ in
any of the expressions $\phi_G(\sigma_{ij})$.  This is a consequence
of the fact that vertex $m$ has no children.  Thus, we have a
partition of the $\sigma_{ij}$ into three sets in which
$\phi_G(\sigma_{ij})$ appear in disjoint sets of variables and there
can be no nontrivial relations involving two or more of these sets of
variables. 
\end{proof}

\begin{prop}
Suppose that for all $i \in [n-1]$, the edge $i \to n \in E(G)$.  Let
$G \setminus n$ be the DAG obtained from $G$ by removing the vertex
$n$.  Then 
$$I_G  =  I_{G \setminus n}  \cdot  \cc[\sigma_{ij} \, \, : \, \, i,j
  \in [n] ].$$
\end{prop}

\begin{proof}
Every vertex in $G \setminus n$ is connected to $n$ and is a parent of
$n$.  This implies that $n$ cannot appear in any conditional
independence statement implied by $G$.  Furthermore, if $C$
$d$-separates $A$ from $B$ in $G \setminus n$, it will $d$-separate
$A$ from $B$ in $G$, because $n$ is below every vertex in $G \setminus
n$.  This implies that the $CI$ statements that hold for $G$ are
precisely the same independence statements that hold for $G \setminus
n$.  Thus 
$$V(C_G) \cap PD_n = V(C_{G \setminus n} \cdot
\cc[\sigma_{ij} \, \, | \, \, i,j \in [n] ]) \cap PD_n.$$
Since $I_G  = I(V(C_G)\cap PD_n)$, this implies the desired equality.
\end{proof}


\section{Tetrad Representation Theorem}  \label{sec:tetrad}

An important step towards understanding the ideals $I_G$ is to derive
interpretations of the polynomials in $I_G$.  We have an
interpretation for a large part of $I_G$, namely, the subideal $C_G
\subseteq I_G$.  Conversely, we can ask when polynomials of a given
form belong to the ideals $I_G$.  Clearly, any linear
polynomial in $I_G$ is a linear combination of polynomials of the form
$\sigma_{ij}$ with $i \neq j$, all of which must also belong to
$I_G$.  Each linear polynomial $\sigma_{ij}$ corresponds to the independence
statement $X_i  \ind X_j$.  Combinatorially, 
the linear from $\sigma_{ij}$ is in $I_G$ if and only if there is no trek
from $i$ to $j$ in $G$.  

A stronger result of this form is the tetrad
representation theorem, first proven in \cite{Spirtes2000}, which
gives a combinatorial characterization of when a tetrad difference 
 $$ \sigma_{ij}  \sigma_{kl}  -  \sigma_{il} \sigma_{jk}$$
belongs to the ideal $I_G$.  The constraints do not necessarily
correspond to conditional independence statements, and need not
belong to the ideal $C_G$.  This will be illustrated in Example \ref{ex:a139}.

The original proof of the tetrad representation theorem in
\cite{Spirtes2000} is quite long and technical.
Our goal in this section is to show how our algebraic perspective can
be used to greatly simplify the proof.  We also include this result
here because we will need the tetrad representation theorem in
Section \ref{sec:tree}. 

\begin{defn}
A vertex $c \in V(G)$ is a \emph{choke point} between sets $I$ and $J$ if
every trek from a point in $I$ to a point in $J$ contains $c$ and
either 
\begin{enumerate}
\item $c$ is on the $I$-side of every trek from $I$ to $J$, or
\item $c$ is on the $J$-side of every trek from $I$ to $J$.
\end{enumerate}
The set of all choke points in $G$ between $I$ and $J$ is denoted $C(I,J)$.
\end{defn}

\begin{ex}
In the graph $c$ is a choke point between $\{1,4\}$ and $\{2,3\}$, but
is not a choke point between $\{1,2\}$ and $\{3,4\}$.
\begin{center}
\resizebox{5cm}{!}{\includegraphics{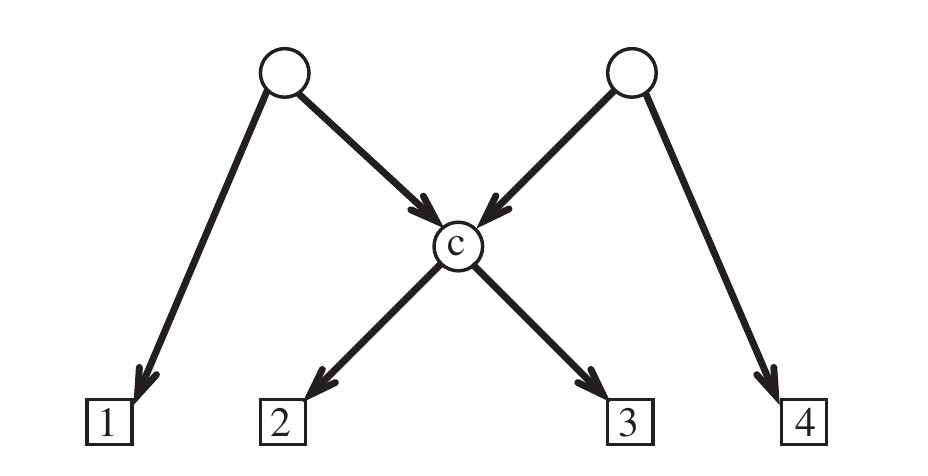}}
\end{center}

\end{ex}

\begin{thm}[Tetrad Representation Theorem \cite{Spirtes2000}]  \label{thm:main}
The tetrad constraint $ \sigma_{ij}  \sigma_{kl}  -  \sigma_{il}
\sigma_{jk}  \, \, = \, \, 0 $ holds for all covariance matrices in
the Bayesian network associated to
$G$ if and only if there is a choke point in $G$ between $\{i,k\}$ and
$\{j,l\}$. 
\end{thm}

Our proof of the tetrad representation theorem will follow after a few lemmas that lead to the irreducible factorization of the polynomials $\sigma_{ij}(a, \lambda)$.

\begin{lemma}
In a fixed DAG $G$, every trek from $I$ to $J$ is incident to every
choke point in $C(I,J)$ and they must be reached always in the same order. 
\end{lemma}
\begin{proof}
If two choke points are on, say, the $I$ side of every trek from $I$
to $J$ and there are two treks which reach these choke points in
different orders, there will be a directed cycle in $G$.  If the choke
points $c_1$ and $c_2$ were on the $I$ side and $J$ side,
respectively, and there were two treks from $I$ to $J$ that reached
them in a different order, this would contradict the property of being
a choke point.
\end{proof}

\begin{lemma}
Let $i = c_0, c_1, \ldots, c_k  = j$ be the ordered choke points in
$C(\{i\}, \{j\})$.   Then the irreducible factorization of
$\sigma_{ij}(a, \lambda)$ is 
$$\sigma_{ij}(a,\lambda) \,  \,  = \, \,  \prod_{t =1}^k   f^t_{ij}(a,
\lambda)$$ 
where $f^t_{ij}(a, \lambda)$ only depends on $\lambda_{pq}$ such that
$p$ and $q$ are between choke points $c_{t-1}$ and $c_{t}$. 
\end{lemma}

\begin{proof}
First of all, we will show that $\sigma_{ij}(a,\lambda)$ has a factorization as
indicated.  Then we will show that the factors  are irreducible.
Define 
$$f^t_{ij}(a, \lambda)  =   \sum_{P \in T(i,j; c_{t-1}, c_{t}) }
a_{{\rm top}(P)}  \prod_{k \to l \in P}  \lambda_{kl} $$ 
where $T(i,j; c_{t-1}, c_{t})$ consists of all paths from
$c_{t-1}$ to $c_{t}$ that are partial treks from $i$ to $j$ (that is,
that can be completed to a trek from $i$ to $j$) and $a_{{\rm top}(P)}
=1$ if the top of the partial trek $P$ is not the top.  When deciding whether or not the top is included in the partial trek, note that almost all choke points are associated with either the $\{i\}$ side or the $\{j\}$ side.  So there is a natural way to decide if $a_{{\rm top}(P)}$ is included or not.  In the exceptional case that $c$ is a choke point on both the $\{i\}$ and the $\{j\}$ side, we repeat this choke point in the list.  This is because $c$ must be the top of every trek from $i$ to $j$, and we will get a factor $f^t_{ij}(a,\lambda) = a_c$.

Since each
$c_{t}$ is a choke point between $i$ and $j$, the product of the
monomials, one from each $f^t_{ij}$, is the monomial corresponding to a
trek from $i$ to $j$.  Conversely, every monomial arises as such a
product in a unique way.  This proves that the desired factorization
holds. 

Now we will show that each of the $f^t_{ij}(a, \lambda)$ cannot
factorize further.  Note that every monomial in $f^t_{ij}(a, \lambda)$
is squarefree in all the $a$ and $\lambda$ indeterminates.  This means
that every monomial appearing in $f^t_{ij}(a, \lambda)$ is a vertex of
the Newton polytope of $f^t_{ij}(a, \lambda)$.  This, in turn, implies
that in any factorization $f^t_{ij}(a, \lambda) = fg$ there is no
cancellation since in any factorization of any polynomial, the
vertices of the Newton polytope is the product of two vertices of the
constituent Newton polytopes.  This means that in any factorization  
$f^t_{ij}(a, \lambda) = fg$, $f$ and $g$ can be chosen to be the sums
of squarefree monomials all with coefficient $1$.   

Now let $f^t_{ij}(a, \lambda) = fg$ be any factorization and let $m$
be a monomial appearing  
in $f^t_{ij}(a, \lambda)$.  If the factorization is nontrivial $m =
m_f m_g$ where $m_f$ and $m_g$ are monomials in $f$ and $g$
respectively.  Since the factorization is nontrivial and $m$
corresponds to a partial trek $P$ in $T(i,j; c_{t-1}, c_t)$, there
must exist a $c$ on $P$ such that, without loss of generality such
that $\lambda_{pc}$ appears in $m_f$ and $\lambda_{cq}$ appears in
$m_g$.  Since every monomial in the expansion of $fg$ corresponds to a
partial trek from $c_{t-1}$ to $c_t$ it must be the case that every
monomial in $f$ contains an indeterminate $\lambda_{sc}$ from some $s$
and similarly, every monomial appearing in $g$ contains a
$\lambda_{cs}$ for some $s$.  But this implies that every partial trek
from $c_{t-1}$ to $c_t$ passes through $c$, with the same
directionality, that is, it is a choke point between $i$ and $j$.
However, this contradicts the fact the $C(\{i\}, \{j\}) = \{c_0,
\ldots, c_t\}$. 
\end{proof}

\noindent {\em Proof of Thm \ref{thm:main}.}   Suppose that the
vanishing tetrad condition holds, that is, $$\sigma_{ij} \sigma_{kl} =
\sigma_{il} \sigma_{kj}$$ 
for all covariance matrices in the model.  This factorization must
thus also hold when we substitute the polynomial expressions in the
parametrization: 
$$\sigma_{ij}(a, \lambda) \sigma_{kl}(a, \lambda) =  \sigma_{il}(a,
\lambda) \sigma_{kj}(a, \lambda).$$ 
Assuming that none of these polynomials are zero (in which case the
choke condition is satisfied for trivial reasons), this means that
each factor $f^t_{ij}(a, \lambda)$ must appear on both the left and
the right-hand sides of this expression.  This is a consequence of the
fact that polynomial rings over fields are unique factorization
domains.  The first factor $f^1_{ij}(a, \lambda)$ could only be a
factor of $\sigma_{il}(a, \lambda)$.  There exists a unique $t \geq 1$
such that $f^1_{ij} \cdots f^{t}_{ij}$ divides $\sigma_{il}$ but
$f^1_{ij} \cdots f^{t+1}_{ij}$ does not divide $\sigma_{il}$.  This
implies that $f^{t+1}_{ij}$ divides $\sigma_{kj}$.  However, this
implies that $c_t$ is a choke point between $i$ and $j$, between $i$
and $l$, between $k$ and $j$.  Furthermore, this will imply that $c_t$
is a choke point between $k$ and $l$ as well, which implies that $c_t$
is a choke point between $\{i,k\}$ and $\{j,l\}$. 

Conversely, suppose that there is a choke point $c$ between $\{i,k\}$
and $\{j,l\}$.  Our unique factorization of the $\sigma_{ij}$ implies
that we can write 
$$\sigma_{ij} =  f_1g_1,  \sigma_{kl}  =  f_2g_2,  \sigma_{il}  =
f_1g_2,  \sigma_{kj}  = f_2g_1$$ 
where $f_1$ and $f_2$ corresponds to partial treks from $i$ to $c$ and
$k$ to $c$, respectively, and $g_1$ and $g_2$ correspond to partial
treks from $c$ to $j$ and $l$, respectively.  Then we have 
$$\sigma_{ij} \sigma_{kl} \, \,  = \, \, f_1g_1f_2g_2  \, \, = \, \,
\sigma_{il} \sigma_{kj},$$ 
so that $\Sigma$ satisfies the tetrad constraint.  \qed

At first glance, it is tempting to suggest that the tetrad
representation theorem says that a tetrad vanishes for every
covariance matrix in the model if and only if an associated
conditional independence statement holds.  Unfortunately, this is not
true, as the following example illustrates.

\begin{ex}\label{ex:a139}
Let $A_{139}$ be the graph with edges $1 \to 4$, $1 \to 5$, $2 \to 4$,
$3 \to 4$ and $4 \to 5$.  Then $4$ is a choke point between $\{2,3\}$
and $\{4,5\}$ and the tetrad $\sigma_{24} \sigma_{35} - \sigma_{25}
\sigma_{34}$ belongs to $I_{A_{139}}$.   However, it is not implied by
the conditional independence statements implied by the graph (that is,  
$\sigma_{24} \sigma_{35} - \sigma_{25}
\sigma_{34} \notin C_{A_{139}}$).  It is precisely this extra tetrad
constraint that forces $A_{139}$ onto the list of graphs that satisfy
$C_G \neq I_G$ from Section \ref{sec:comp}.
\end{ex}

In particular, a choke point between two sets need not be a
$d$-separator of those sets.  In the case
that $G$ is a tree, it is true 
that tetrad constraints are conditional independence constraints.

\begin{prop}\label{prop:treeind}
Let $T$ be a tree and suppose that $c$ is a choke point between $I$
and $J$ in $T$.  Then either $c$ $d$-separates $I \setminus \{c\}$ and
$J \setminus \{c\}$ or $\emptyset$ $d$-separates $I \setminus \{c\}$
and $J \setminus \{c\}$.
\end{prop}

\begin{proof}
Since $T$ is a tree, there is a unique path from an element in $I
\setminus c$ to an element in $J \setminus c$.  If this path is not a
trek, we have $\emptyset$ $d$-separates $I \setminus \{c\}$ from $J
\setminus \{c\}$.  On the other hand, if this path is always a trek we
see that $\{c\}$ $d$-separates $I \setminus \{c\}$ from $J
\setminus \{c\}$.
\end{proof}

The tetrad representation theorem gives a simple combinatorial rule
for determining when a $2 \times 2$ minor of $\Sigma$ is in $I_G$.
More generally, we believe that there should exist a graph theoretic 
rule that determines when a general determinant $| \Sigma_{A,B} |
\in I_G$ in terms of structural features of the DAG $G$.  The
technique we have used above, which relies on giving a factorization
of the polynomials $\sigma(a, \lambda)$, does not seem like it will
extend to higher order minors.  One approach at a generalization of
the tetrad representation theorem would be to find a
cancellation free expression for the determinant $|\Sigma_{A,B}|$ in
terms of the parameters $a_i$ and $\lambda_{ij}$, along the lines of
the Gessel-Viennot theorem \cite{Gessel1985}.  From such a result, one
could deduce a 
combinatorial rule for when $|\Sigma_{A,B}|$ is zero.  This suggests the
following problem. 

\begin{pr}
Develop a Gessel-Viennot theorem for treks; that is, determine a
combinatorial formula for the expansion of $|\Sigma_{A,B}|$ in terms
of the treks in $G$.
\end{pr}


\section{Fully Observed Trees}  \label{sec:tree}

In this section we study the Bayesian networks of trees in the
situation where all random variables are observed.  We show that the
toric ideal $I_T$ is generated by linear forms $\sigma_{ij}$ and
quadratic tetrad constraints.  The Tetrad Representation Theorem and
Proposition \ref{prop:treeind} then
imply that $I_T = C_T$.   We also investigate further algebraic
properties of the ideals $I_T$ using the fact that $I_T$ is a toric
ideal and some techniques from polyhedral geometry.

For the rest of this section, we assume that $T$ is a tree, where by a
tree we mean a DAG whose underlying undirected graph is a tree.  These
graphs are sometimes called polytrees in the graphical models
literature.  A \emph{directed tree} is a tree all of whose edges are
directed away from a given source vertex. 

Since $I_T$ is a toric ideal, it can be analyzed using techniques from
polyhedral geometry.  In particular, for each $i,j$ such that $T(i,j)$
is nonempty, let $\bfa_{ij}$ denote the exponent vector of the
monomial $\sigma_{ij} = a_{\rmtop(P)} \prod_{k \to l \in P}
\lambda_{kl}$.  Let $A_T$ denote the set of all these exponent
vectors.  The geometry of the toric variety $V(I_T)$ is determined by
the discrete geometry of the polytope $P_T =  \conv(A_T)$. 

The polytope $P_T$ is naturally embedded in $\rr^{2n-1}$, where $n$ of
the coordinates on 
$\rr^{2n-1}$ correspond to the vertices of $T$ and $n-1$ of the
coordinates correspond to the edges of $T$.  Denote the first set of
coordinates by $x_i$ and the second by $y_{ij}$ where $i \to j$ is an
edge in $T$.  Our first results is a description of the facet
structure of the polytope $P_T$. 

\begin{thm}\label{thm:facets}
The polytope $P_T$ is the solution to the following set of equations
and inequalities: 

\smallskip

\begin{tabular}{ll}
$x_i \geq 0$  &  for all $i \in V(T)$ \\
$y_{ij} \geq 0$ &  for all $i \to j \in E(T)$  \\
$\sum_{i  \in V(T)} x_i =1$  &   \\
$x_j + \sum_{i:  \, \, i \to j \in E(T)} y_{ij}  \,\,  - \,\, y_{jk} \, \, \geq \,\, 0 $   &  for all $j \to k \in E(T)$ \\
$2x_j + \sum_{i:  \, \, i \to j \in E(T)}  y_{ij} \, \, - \,\,   \sum_{k: \, \, j \to k \in E(T)}  y_{jk} \, \, \geq \, \, 0 $  &  for all $j \in V(T).$ 
\end{tabular}
\end{thm}

\begin{proof}
Let $Q_T$ denote the polyhedron defined as the solution space to the
given constraints.   First of all, $Q_T$ is bounded.  To see this, first note that because of the positive constraints and the equation $\sum_{i  \in V(T)} x_i =1$, we have that $x_i \leq 1$ is implied by the given constraints.   Then, starting from the sources of the tree and working our way down the edges repeatedly using the inequalities $x_j + \sum_{i:  \, \, i \to j \in E(T)} y_{ij}  \,\,  - \,\, y_{jk} \, \, \geq \,\, 0 $, we see that the $y_{ij}$ are also bounded.

Now, we have $P_T \subseteq Q_T$, since
every trek will satisfy any of the indicated constraints.  Thus, we
must show that $Q_T \subseteq P_T$.  To do this, it suffices to show
that for any vector $(x^0, y^0) \in Q_T$, there exists $\lambda >0$,
$(x^1, y^1)$ and $(x^2, y^2)$ such that  
$$(x^0,y^0) = \lambda(x^1,y^1) + (1 - \lambda)(x^2, y^2)$$
where $(x^1, y^1)$ is one of the $0/1$ vectors $\bfa_{ij}$
and $(x_2, y_2) \in Q_T$.  Because $Q_T$ is bounded, this  will imply that the extreme points of $Q_T$ are a subset of the extreme points of  $P_T$, and hence $Q_T \subseteq P_T$.   Without loss of generality we may suppose that all of the 
coordinates $y^0_{ij}$ are positive, otherwise the problem reduces to a
smaller tree or forest because the resulting inequalities that arise
when $y_{ij} =0$ are precisely those that are necessary for the
smaller tree.  Note that for a forest $F$, the polytope $P_F$ is the
direct join of polytopes $P_T$ as $T$ ranges over the connected
components of $F$, by Proposition \ref{prop:disconnect}. 

For any fixed $j$, there cannot exist distinct values $k_1$, $k_2$,
and $k_3$ such that all of 
$$x^0_j + \sum_{i:\, \, i \to j \in E(T)} y^0_{ij} -  y^0_{jk_1}  = 0 $$
$$x^0_j + \sum_{i:\, \, i \to j \in E(T)} y^0_{ij}  -  y^0_{jk_2}  = 0 $$
$$x^0_j + \sum_{i:\, \, i \to j \in E(T)} y^0_{ij} -  y^0_{jk_3}  = 0 $$
hold.  If there were, we could add these three equations together to
deduce that 
$$3x^0_j + 3 \sum_{i:\, \, i \to j \in E(T)} y^0_{ij}  -  y^0_{jk_1} -
y^0_{jk_2} - y^0_{jk_3} = 0.$$ 
This in turn implies that
$$2x^0_j +  \sum_{i:\, \, i \to j \in E(T)} y^0_{ij}  -  y^0_{jk_1} -
y^0_{jk_2} - y^0_{jk_3} \leq 0$$ 
with equality if and only if $\rmpa(j) = \emptyset$ and $x^0_j = 0$.
This in turn implies that, for instance,  $y^0_{jk_1} =0$
contradicting our assumption that $y^0_{ij} > 0$ for all $i$ and $j$.
By a similar argument, if exactly two of these facet defining
inequalities hold sharply, we see that 
$$2x^0_j + 2 \sum_{i:\, \, i \to j \in E(T)} y^0_{ij}  -  y^0_{jk_1} -
y^0_{jk_2}  = 0$$ 
which implies that $j$ has exactly two descendants and no parents.

Now mark each edge $j \to k$  in the tree $T$ such that
$$x^0_j + \sum_{i:\, \, i \to j \in E(T)} y^0_{ij} -  y^0_{jk}  =  0. $$
By the preceding paragraph, we can find a trek $P$ from a sink in the
tree to a source in the tree and (possibly) back to a different sink
that has the property that for no $i$ in the trek there exists $k$ not
in the path such that $i \to k$ is a marked edge.  That is, the
preceding paragraph shows that there can be at most $2$ marked edges
incident to any given vertex. 

Given $P$, let $(x^1, y^1)$ denote the corresponding $0/1$ vector.  We
claim that there is a $\lambda >0$ such that  
\begin{equation}\label{eq:define}
(x^0,y^0) = \lambda(x^1,y^1) + (1 - \lambda)(x^2, y^2)
\end{equation}
holds with $(x^2, y^2) \in Q_T$.  Take $\lambda>0$ to be any very
small number and define $(x^2, y^2)$ by the given equation.  Note that
by construction the inequalities $x^2_i \geq 0$ and $y^2_{ij} \geq 0$
will be satisfied since for all the nonzero entries in $(x^1, y^1)$,
the corresponding inequality for $(x^0, y^0)$ must have been nonstrict
and $\lambda$ is small.  Furthermore, the constraint $\sum x^2_i =1$
is also automatically satisfied.  It is also easy to see that the
last set of inequalities will also be satisfied since through each
vertex the path will either have no edges, an incoming edge and an
outgoing edge, or two outgoing edges and the top vertex, all of which
do not change the value of the linear functional. 

Finally to see that the inequalities of the form 
$$x_j + \sum_{i:  \, \, i \to j \in E(T)} y_{ij}  \,\,  - \,\, y_{jk}
\, \, \geq \,\, 0 $$ 
are still satisfied by $(x^2,y^2)$, note that marked edges of $T$ are
either contained in the path $P$ or not incident to the path $P$.
Thus, the strict inequalities remain strict (since they will involve
modifying by an incoming edge and an outgoing edge or an outgoing edge
and the top vertex), and the nonstrict inequalities remain nonstrict
since $\lambda$ is small.  Thus, we conclude that $Q_T \subseteq P_T$,
which completes the proof.
\end{proof}

\begin{cor}\label{cor:pull}
Let $\prec$ be any reverse lexicographic term order such that
$\sigma_{ii}  \succ \sigma_{jk}$ for all $i$ and $ j \neq k$.  Then
${\rm in}_\prec(I_T)$ is squarefree.  In other words, the associated
pulling triangulation of $P_T$ is unimodular.  
\end{cor}

\begin{proof}
The proof is purely polyhedral, and relies on the geometric
connections between triangulations and initial ideals of toric ideals.
See Chapter 8 in \cite{Sturmfels1996} for background on this material
including pulling triangulations.  Let $\bfa_{ij}$ denote the vertex
of $P_T$ corresponding to the monomial $\phi_G(\sigma_{ij})$.  For $i
\neq j$, each of the vertices $\bfa_{ij}$ has lattice distance at most
one from any of the facets described by Theorem \ref{thm:facets}.
This is seen by evaluating each of the linear functionals at the $0/1$
vector corresponding to the trek between $i$ and $j$.   

If we pull from one of these vertices we get a unimodular
triangulation provided that the induced pulling triangulation on each
of the facets of $P_T$ not containing $\bfa_{ij}$ is unimodular.  This
is because the normalized volume of a simplex is the volume of the
base times the lattice distance from the base to the vertex not on the
base. 

The facet defining inequalities of any face of $P_T$ are obtained by
taking an appropriate subset of the facet defining inequalities of
$P_T$.  Thus, as we continue the pulling triangulation, if the current
face contains a vertex $\bfa_{ij}$ with $i \neq j$, we will pull from
this vertex first and get a unimodular pulling triangulation provided
the induced pulling triangulation of every face is unimodular.  Thus,
by induction, it suffices to show that the faces of $P_T$ that are the
convex hull of vertices $\bfa_{ii}$ have unimodular pulling
triangulations.  However, these faces are always unimodular
simplices. 
\end{proof}

\begin{cor}
The ring $\cc[\Sigma]/I_T$ is normal and Cohen-Macaulay when $T$ is a tree.
\end{cor}

\begin{proof}
Since $P_T$ has a unimodular triangulation, it is a normal polytope
and hence the semigroup ring $\cc[\Sigma]/I_T$ is normal.  Hochster's
theorem \cite{Hochster1972} then implies that $\cc[\Sigma]/I_T$ is
Cohen-Macaulay.  
\end{proof}

While we know that $\cc[\Sigma]/I_T$ is always Cohen-Macaulay, it
remains to determine how the Cohen-Macaulay type of $I_T$ depends on
the underlying tree $T$.  Here is a concrete conjecture concerning
the special case of Gorenstein trees. 

\begin{conj}
Suppose that $T$ is a directed tree.  Then $\cc[\Sigma]/I_T$ is
Gorenstein if and only if the degree of every vertex in $T$ is less
than or equal to three. 
\end{conj}

A \emph{downward directed tree} is a tree all of whose edges point to the
unique sink in the tree.  A leaf of such a downward directed tree is
then a source of the tree.  With a little more refined information
about which inequalities defining $P_T$ are facet defining, we can
deduce results about the degrees of the ideals $I_T$ in some cases. 

\begin{cor}\label{cor:treedegree}
Let $T$ be a downward directed tree and let $i$ be any leaf of $T$,
$s$ the sink of $T$, and $P$ the unique trek in $T(i,s)$.  Then 
$$\deg I_T  \quad = \quad  \sum_{k \to l \in P}  \deg I_{T \setminus
  k \to l}$$ 
where $T \setminus k \to l$ denotes the forest obtained from $T$
by removing the edge $k \to l$. 
\end{cor}

\begin{proof}
First of all, note that in the case of a downward directed tree the
inequalities of the form 
$$2x_j + \sum_{i:  \, \, i \to j \in E(T)}  y_{ij} \, \, - \,\,
\sum_{k: \, \, j \to k \in E(T)}  y_{jk} \, \, \geq \, \, 0 $$ 
are redundant:  since each vertex has at most one descendant, it is
implied by the the other constraints. 
Also, for any source $t$, the inequality $x_t \geq 0$ is redundant,
because it is implied by the inequalities $x_t -  y_{tj}  \geq 0$ and
$y_{tj} \geq 0$ where $j$ is the unique child of $t$.   

Now we will compute the normalized volume of the polytope $P_T$ (which
is equal to the degree of the toric ideal $I_T$) by computing the
pulling triangulation from Corollary  \label{cor:pull} and relating
the volumes of the pieces to the associated subforests. 

Since the pulling triangulation of $P_T$ with $\bfa_{is}$ pulled first
is unimodular, the volume of $P_T$ is the sum of the volumes of the
facets of $P_T$ that do not contain $\bfa_{is}$.  Note that
$\bfa_{is}$ lies on all the facets of the form 
$$x_j + \sum_{i:  \, \, i \to j \in E(T)} y_{ij}  \,\,  - \,\, y_{jk}
\, \, \geq \,\, 0 $$ 
since through every vertex besides the source and sink, the trek has
either zero or
two edges incident to it.  Thus, the only facets that $\bfa_{is}$ does
not lie on are of the form $y_{kl} \geq 0$ such that $k \to l$ is an
edge in the trek $P$.  However, the facet of $P_T$ obtained by setting
$y_{kl} = 0$ is precisely the polytope $P_{T \setminus k \to l}$,
which follows from Theorem
\ref{thm:facets}. 
\end{proof}

Note that upon removing an edge in a tree we obtain a forest.
Proposition \ref{prop:disconnect} implies that the degree of such a
forest is the product of the degrees of the associated trees.  Since
the degree of the tree consisting of a single point is one, the
formula from Corollary \ref{cor:treedegree} yields a recursive
expression for the degree of a downward directed forest. 

\begin{cor}
Let $T_n$ be the directed chain with $n$ vertices.  Then $\deg I_{T_n}
=  \frac{1}{n}{2n-2 \choose n-1}$, the $n-1$st Catalan number.  
\end{cor}

\begin{proof}
In Corollary \ref{cor:treedegree} we take the unique path from $1$ to
$n$.  The resulting forests obtained by removing an edge are the
disjoint unions of two paths.  By the product formula implied by
Proposition \ref{prop:disconnect} we deduce that the degree of
$I_{T_n}$ satisfies the recurrence: 
$$ \deg I_{T_n}    =  \sum_{i =1}^{n-1}   \deg I_{T_i}  \cdot  \deg
I_{T_{n-i}}$$ 
with initial condition $\deg I_{T_1} = 1$.  This is precisely the
recurrence and initial conditions for the Catalan numbers \cite{Stanley1999}. 
\end{proof}

Now we want to prove the main result of this section, that the
determinantal conditional independence statements actually generate
the ideal $I_T$ when $T$ is a tree.  To do this, we will exploit the
underlying toric structure, introduce a tableau notation for working
with monomials, and introduce an appropriate ordering of the
variables.  

Each variable $\sigma_{ij}$ that is not zero can be identified with
the unique trek in $T$ from $i$ to $j$.  We associate to
$\sigma_{ij}$ the tableau which records the elements of $T$ in this unique
trek, which is represented like this: 
$$\sigma_{ij}  =  [\underline{a}Bi |  \underline{a}Cj]$$
where $B$ and $C$ are (possibly empty) strings.  If, say, $i$ were at
the top of the path, we would write the tableau as 
$$\sigma_{ij} = [\underline{i} |  \underline{i}Cj].$$
The tableau is in its standard form if $aBi$ is lexicographically
earlier than $aCj$.  We introduce a lexicographic total order on
standard  form tableau variables by declaring 
$[\underline{a}A|\underline{a}B] \prec
[\underline{c}C|\underline{c}D]$ if $aA$ is lexicographically smaller
that $cC$, or if $aA =cC$ and $aB$ is lexicographically smaller than
$cD$.  Given a monomial, its tableau representation is the row-wise
concatenation of the tableau forms of each of the variables appearing
in the monomial.   

\begin{ex}
Let $T$ be the tree with edges $1 \to 3$, $1 \to 4$, $2 \to 4$, $3 \to
5$, $3\to 6$, $4 \to 7$, and $4 \to 8$.  Then the monomial
$\sigma_{14}\sigma_{18}\sigma_{24}\sigma^2_{34}\sigma_{38}
\sigma_{57}\sigma_{78}$    
has the standard form lexicographically ordered tableau:
$$\left[\begin{array}{l|l}
\underline{1}  &  \underline{1}4 \\
\underline{1} & \underline{1}48 \\
\underline{1}3 & \underline{1}4 \\
\underline{1}3 & \underline{1}4 \\
\underline{1}3 & \underline{1}48 \\
\underline{1}35 & \underline{1}47 \\
\underline{2} & \underline{2} 4 \\
\underline{4}7 & \underline{4}8 \end{array}  \right].$$
Note that if a variable appears to the $d$-th power in a monomial, the
representation for this variable is repeated as $d$ rows in the
tableau.  \qed 
\end{ex}

When we write out general tableau, lower-case letters will always
correspond to single characters (possibly empty) and upper case
letters will always correspond to strings of characters (also,
possibly empty).

\begin{thm}
For any tree $T$, the conditional independence statements implied by
$T$ generate $I_T$.  In particular, $I_T$ is generated by linear
polynomials $\sigma_{ij}$ and quadratic tetrad constraints. 
\end{thm}

\begin{proof}
First of all, we can ignore the linear polynomials as they always
correspond to independence constraints and work modulo these linear
constraints when working with the toric ideal $I_T$.  In addition,
every quadratic binomial of the form $\sigma_{ij}\sigma_{kl} - 
\sigma_{il}\sigma_{kj}$ that belongs to $I_T$ is implied by a
conditional independence statement.  This follows from Proposition
\ref{prop:treeind}.  
Note that this holds even if the set $\{i,j,k,l\}$ does not have four elements.
Thus, it suffices to show that $I_T$ modulo the linear constraints is
generated by quadratic binomials. 

To show that $I_T$ is generated by quadratic binomials, it suffices to
show that any binomial in $I_T$ can be written as a polynomial linear
combination of the quadratic binomials in $I_T$.  This, in turn, will
be achieved by showing that we can ``move'' from the tableau
representation of one of the monomials to the other by making local
changes that correspond to quadratic binomials.  To show this last
part, we will define  a sort of distance between two monomials and show
that it is always possible to decrease this distance using these
quadratic binomials/ moves.  This is a typical trick for dealing with
toric ideals, illustrated, for instance, in \cite{Sturmfels1996}.

To this end let $f$ be a binomial in $I_T$.  Without loss of
generality, we may 
suppose the terms of $f$ have no common factors, because if
$\sigma^\bfa \cdot f  \in I_T$ then $f \in I_T$ as well.  We will
write $f$ as the difference of two tableaux, which are in standard
form with their rows lexicographically ordered.  The first row in the
two tableaux are different and they have a left-most place where they
disagree.  We will show that we can always move this position further
to the right.  Eventually the top rows of the tableaux will agree and
we can delete this row (corresponding to the same variable) and arrive
at a polynomial of smaller degree. 

Since $f \in I_T$, the treks associated to the top rows of the two
tableaux must 
have the same top.  There are two cases to consider.  Either the first
disagreement is immediately after the top or not.  In the first case,
this means that the binomial $f$ must have the form: 
$$\left[ 
\begin{array}{l|l}
\underline{a}bB  &  \underline{a}cC  \\
&   \end{array} \right]
-
\left[
\begin{array}{l|l}
\underline{a}bB  &  \underline{a}dD \\
&  \end{array} \right].
$$
Without loss of generality we may suppose that $c<d$.  Since $f \in
I_T$ the string $ac$ must appear somewhere on the right-hand monomial.
Thus, $f$ must have the form: 
$$\left[ 
\begin{array}{l|l}
\underline{a}bB  &  \underline{a}cC  \\
 &   \\
  &  \\
\end{array} \right]
-
\left[
\begin{array}{l|l}
\underline{a}bB  &  \underline{a}dD \\
\underline{a}eE  &  \underline{a}cC' \\
&  \end{array} \right].
$$
If $d \neq e$, we can apply the quadratic binomial
$$\left[
\begin{array}{l|l}
\underline{a}bB  &  \underline{a}dD \\
\underline{a}eE  &  \underline{a}cC'  
\end{array} \right] -
\left[
\begin{array}{l|l}
\underline{a}bB  &  \underline{a}cC' \\
\underline{a}eE  &  \underline{a}dD \\
\end{array} \right]
$$
to the second monomial to arrive at a monomial which has fewer
disagreements with the left-hand tableau in the first row.  On the
other hand, if $d = e$, we cannot apply this move (its application
results in ``variables'' that do not belong to $\cc[\Sigma]$).  Keeping
track of all the $ad$ patterns that appear on the right-hand side, and
the consequent $ad$ patterns that appear on the left-hand side, we
see that our binomial $f$ has the form 

$$\left[ 
\begin{array}{l|l}
\underline{a}bB  &  \underline{a}cC  \\
\underline{a}d* & *  \\
\vdots   &  \vdots  \\ 
\underline{a}d* &  *
\end{array} \right]
-
\left[
\begin{array}{l|l}
\underline{a}bB  &  \underline{a}dD \\
\underline{a}dD'  &  \underline{a}cC' \\
\underline{a}d*   &   *   \\
\vdots  &  \vdots  \\
\underline{a}d*  &  * 
  \end{array} \right].
$$ 
Since there are the same number of $ad$'s on both sides we see that
there is at least one more $a$ on the right-hand side which has no
$d$'s attached to it.  Thus, omitting the excess $ad$'s on both sides,
our binomial $f$ contains: 
$$\left[ 
\begin{array}{l|l}
\underline{a}bB  &  \underline{a}cC  \\
 &   \\
   &    \\   
\end{array} \right]
-
\left[
\begin{array}{l|l}
\underline{a}bB  &  \underline{a}dD \\
\underline{a}dD'  &  \underline{a}cC' \\
\underline{a}eE &   \underline{a}gG   \\
\end{array} \right].
$$ 
with $d \neq e$ or $g$.  We can also assume that $c\neq e,g$
otherwise, we could apply a quadratic move as above. 
Thus we apply the quadratic binomials
$$\left[
\begin{array}{l|l}
\underline{a}dD'  &  \underline{a}cC' \\
\underline{a}eE  &  \underline{a}gG  
\end{array} \right] -
\left[
\begin{array}{l|l}
\underline{a}dD'  &  \underline{a}gG\\
\underline{a}eE  &  \underline{a}cC'\\
\end{array} \right]
$$
and
$$\left[
\begin{array}{l|l}
\underline{a}bB  &  \underline{a}dD \\
\underline{a}eE  &  \underline{a}cC'  
\end{array} \right] -
\left[
\begin{array}{l|l}
\underline{a}bB  &  \underline{a}cC' \\
\underline{a}eE  &  \underline{a}dD \\
\end{array} \right]
$$
to reduce the number of disagreements in the first row.  This
concludes the proof of the first case.  Now suppose that the first
disagreement does not occur immediately after the $a$.  Thus we may
suppose that $f$ has the form: 

$$
\left[
\begin{array}{l|l}
\underline{a}AxbB &  \underline{a}C  \\
 &    
 \end{array}  \right]  -
 \left[
 \begin{array}{l|l}
 \underline{a}AxdD &  \underline{a}E  \\
  &  
  \end{array}  \right].
$$
Note that it does not matter whether or not this disagreement appears
on the left-hand or right-hand side of the tableaux.  Since the string
$xd$ appears on right-hand monomial it must also appear somewhere on
the left-hand monomial as well.  If $x$ is not the top in this
occurrence, we can immediately apply a quadratic binomial to reduce the
discrepancies in the first row.  So we may assume the $f$ has the
form: 
$$
\left[
\begin{array}{l|l}
\underline{a}AxbB &  \underline{a}C  \\
\underline{x}dD'  & \underline{x}gG  \\
 &      
 \end{array}  \right]  -
 \left[
 \begin{array}{l|l}
 \underline{a}AxdD &  \underline{a}E  \\
  &    \\
  &   \\
  \end{array}  \right].
$$
If $b \neq g$ we can apply the quadratic binomial
$$\left[
\begin{array}{l|l}
\underline{a}AxbB &  \underline{a}C  \\
\underline{x}dD'  & \underline{x}gG  \\
\end{array}  \right]  -
\left[
\begin{array}{l|l}
\underline{a}AxdD' &  \underline{a}C  \\
\underline{x}bB  & \underline{x}gG  \\
\end{array}  \right]$$
to the left-hand monomial to reduce the discrepancies in the first
row.  So suppose that $g = b$.  Enumerating the $xb$ pairs that can
arise on the left and right hand monomials, we deduce, akin to our
argument in the first case above, that $f$ has the form: 
$$
\left[
\begin{array}{l|l}
\underline{a}AxbB &  \underline{a}C  \\
\underline{x}dD'  & \underline{x}bG  \\
\underline{x}hH & \underline{x}kK  \\ 
&  \\     
 \end{array}  \right]  -
 \left[
 \begin{array}{l|l}
 \underline{a}AxdD &  \underline{a}E  \\
  &    \\
  &   \\
  &  \\
  \end{array}  \right]
$$
where $h$ and $k$ are not equal to $b$ or $d$.  Then we can apply the
two quadratic binomials: 
$$
\left[
\begin{array}{l|l}
\underline{x}dD' & \underline{x}bG \\
\underline{x}hH  &  \underline{x}kK \\
 \end{array}  \right]  -
 \left[
 \begin{array}{l|l}
\underline{x}hH   & \underline{x}bG    \\
 \underline{x}dD' & \underline{x}kK   \\
  \end{array}  \right]
$$
and
$$
\left[
\begin{array}{l|l}
\underline{a}AxbB &  \underline{a}C  \\
\underline{x}dD'  & \underline{x}kK  \\
     
 \end{array}  \right]  -
 \left[
 \begin{array}{l|l}
 \underline{a}AxdD' &  \underline{a}C  \\
\underline{x}bB  & \underline{x}kK    \\
  \end{array}  \right]
$$
to the left-hand monomial to produce a monomial with fewer
discrepancies in the first row.  We have shown that no matter what
type of discrepancy that can occur in the first row, we can always
apply quadratic moves to produce fewer discrepancies.  This implies
that $I_T$ is generated by quadrics.  
\end{proof}

Among the results in this section were our proofs that $I_T$ has a
squarefree initial ideal (and hence $\cc[\Sigma]/I_T$ is normal and
Cohen-Macaulay) and that $I_T$ is generated by linear forms and
quadrics.  It seems natural to wonder if there is a term order that
realizes these two features simultaneously.

\begin{conj}
There exists a term order $\prec$ such that ${\rm in}_\prec(I_T)$ is
generated by squarefree monomials of degree one and two.
\end{conj}


\section{Hidden Trees}  \label{sec:hidtree}

This section and the next concern Bayesian networks with hidden
variables.  A hidden or latent random variable is one which we do not
have direct access to.  These hidden variables might represent
theoretical quantities that are directly unmeasurable (e.g. a random
variable representing intelligence), variables we cannot have access
to (e.g. information about extinct species), or variables that have
been censored (e.g. for sensitive random variables in census data).
If we are given a model over all the observed and hidden random
variables, the partially observed model is the one obtained by
marginalizing over the hidden random variables.  A number of
interesting varieties arise in this hidden variable setting. 

For Gaussian random variables, the marginalization is again Gaussian,
and the mean and covariance matrix are obtained by extracting the
subvector and submatrix of the mean and covariance matrix
corresponding to the observed random variables.   This immediately
yields the following proposition.

\begin{prop}\label{prop:elim}
Let $I  \subseteq \cc[\mu, \Sigma]$ be the vanishing ideal for a
Gaussian model.  Let  $H \cup O =[n]$ be a partition of the
random variables into hidden and observed variables $H$ and $O$.  Then 
$$I_{O}  \quad := \quad  I \cap \cc[ \mu_i, \sigma_{ij} \, \, | \, \,
  i,j \in O]$$ 
is the vanishing ideal for the partially observed model.
\end{prop}

\begin{proof}
Marginalization in the Gaussian case corresponds to projection onto
the subspace of pairs $(\mu_O, \Sigma_{O,O})  \subseteq  \rr^{|O|}
\times  \rr^{ {|O| +1 \choose 2 }}$.  Coordinate projection 
is equivalent to elimination  \cite{Cox1997}. 
\end{proof}

In the case of a Gaussian Bayesian network, Proposition
\ref{prop:elim} has a number of useful corollaries, of both a
computational and theoretical nature.  First of all, it allows for the
computation of the ideals defining a hidden variable model as an easy
elimination step.  Secondly, it can be used to explain the phenomenon
we observed in Example \ref{ex:verma}, that the constraints defining a
hidden variable model appeared as generators of the ideal of the fully
observed model. 

\begin{defn}
Let $H \cup O$ be a partition of the nodes of the DAG $G$.  The hidden
nodes $H$ are said to be \emph{upstream} from the observed nodes $O$
in $G$ if there are no edges $o \to h$ in $G$ with $o \in O$ and $h \in
H$. 
\end{defn}

If $H \cup O$ is an upstream partition of the nodes of $G$, we
introduce a grading on the ring $\cc[a, \lambda]$ which will, in
turn, induce a grading on $\cc[\Sigma]$. 
Let $\deg a_h = (1,0)$ for all $h \in H$, $\deg a_o = (1,2)$ for all
$o \in O$,  $\deg \lambda_{ho}  =  (0,1)$ if $h \in H$ and $o \in O$,
and $\deg \lambda_{ij} = (0,0)$ otherwise. 

\begin{lemma}
Suppose that $H \cup O = [n]$ is an upstream partition of the vertices
of $G$.  Then each of the polynomials $\phi_G(\sigma_{ij})$ is
homogeneous with respect to the upstream grading and  
$$\deg(\sigma_{ij}) =  \left\{  \begin{array}{cl} 
(1,0) & \mbox{ if } i \in H, j \in H \\
(1,1) & \mbox{ if } i \in H, j \in O  \mbox{ or }  i \in O, j \in H\\
(1,2) & \mbox{ if } i \in O, j \in O. \end{array}  \right.$$
Thus, $I_G$ is homogeneous with respect to the induced grading on
$\cc[\Sigma]$. 
\end{lemma} 

\begin{proof}
There are three cases to consider.  If both $i, j \in H$, then every
trek in $T(i,j)$ has a top element in $H$ and no edges of the form $h
\to o$.  In this case, the degree of each path is the vector $(1,0)$.
If $i \in H$ and $j \in O$, every trek from $i$ to $j$
has a top in $H$ and exactly one edge of the form $h \to o$.  Thus,
the degree of every monomial in $\phi(\sigma_{ij}) $ is $(1,1)$.  If
both $i,j \in O$, then either each trek $P$ from $i$ to
$j$ has a top in $O$, or has a top in $H$.  In the first case there
can be no edges in $P$ of the form $h \to o$, and in the second case
there must be exactly two edges in $P$ of the form $h \to o$.  In
either case, the degree of the monomial corresponding to $P$ is
$(1,2)$. 
\end{proof}

Note that the two dimensional grading we have described can be
extended to an $n$ dimensional grading on the ring $\cc[\Sigma]$ by
considering all collections of upstream variables in $G$ simultaneously.

\begin{thm}[Upstream Variables Theorem]\label{thm:upstream}
Let $H \cup O$ be an upstream partition of the vertices of $G$.  Then
every minimal generating set of $I_G$ that is
homogeneous with respect to the upstream grading contains a minimal
generating set of $I_{G,O}$.  
\end{thm}

\begin{proof}
The set of indeterminates $\sigma_{ij}$ corresponding to the observed
variables are precisely the variables whose degrees lie on the facet
of the degree semigroup generated by the vector $(1,2)$.  This implies
that the subring generated by these indeterminates is a facial
subring. 
\end{proof}

The upstream variables theorem is significant because any natural
generating set of an ideal $I$ is homogeneous with respect to its
largest homogeneous grading group.  For instance, every reduced
Gr\"obner basis if $I_G$ will be homogeneous with respect to the
upstream grading.  For trees, the upstream variables theorem immediately implies:

\begin{cor}\label{thm:hidtree}
Let $T$ be a rooted directed tree and $O$ consist of the leaves of
$T$.  Then $I_{T,O}$ is generated by the quadratic tetrad constraints
$$\sigma_{ik} \sigma_{jl} - \sigma_{il} \sigma_{kj}$$
such that $i,j,k,l \in O$, and there is a choke point $c$ between
$\{i,j\}$ and $\{k,l\}$. 
\end{cor}

Corollary \ref{thm:hidtree} says that the ideal of a hidden tree model
is generated by the tetrad constraints induced by the choke points in
the tree.  Sprites et al  \cite{Spirtes2000} use these tetrad
constraints as a tool for inferring DAG models with hidden variables.
Given a sample 
covariance matrix, they test whether a collection of tetrad
constraints is equal to zero.  From the given tetrad constraints that
are satisfied, together with the tetrad representation theorem, they
construct a DAG that is consistent with these vanishing tetrads.
However, it is not clear from that work whether or not it is enough to
consider only these tetrad constraints.  Indeed, as shown in
\cite{Spirtes2000}, there are pairs of graphs with hidden nodes that
have precisely the same set of tetrad constraints that do not yield
the same family of covariance matrices.  Theorem
\ref{thm:hidtree} can be seen as a mathematical justification of the
tetrad procedure of Spirtes, et al, in the case of hidden tree models,
because it shows 
that the tetrad constraints are enough to distinguish between the
covariance matrices coming from different trees.


\section{Connections to Algebraic Geometry}\label{sec:alggeom}

In this section, we give families of examples to show how classical
varieties from algebraic geometry arise in the study of Gaussian
Bayesian networks.  In particular, we show how toric degenerations of
the Grassmannian, matrix Schubert varieties, and secant varieties all
arise as special cases of Gaussian Bayesian networks with hidden
variables.


\subsection{Toric Initial Ideals of the Grassmannian}

Let $Gr_{2,n}$ be the Grassmannian of $2$-planes in $\cc^n$.  The
Grassmannian has the natural structure of an algebraic variety under
the Pl\"ucker embedding.  The ideal of the Grassmannian is generated
by the quadratic Pl\"ucker relations: 
$$I_{2,n} :=  I( Gr_{2,n})  =  \left<  \sigma_{ij} \sigma_{kl}  -
\sigma_{ik} \sigma_{jl} +  \sigma_{il} \sigma_{jk}  \, \, | \, \,1
\leq i< j < k < l \leq n \right> \subset \cc[\Sigma].$$ 

The binomial initial ideals of $I_{2,n}$ are in bijection with the unrooted
trivalent trees with $n$ leaves.   
These binomial initial ideals are, in fact, toric ideals, and we will
show that: 

\begin{thm}\label{thm:grin}
Let $T$ be a rooted directed binary tree with $[n]$ leaves and let $O$ be the set
of leaves of $T$.  Then there is a weight vector
$\omega \in   \rr^{ {n+1 \choose 2}}$ and a sign vector $\tau \in
\{\pm 1 \}^{ { n+1 \choose 2}}$ such that 
$$ I_{T,O}  =   \tau \cdot  {\rm in}_\omega (I_{2,n}).$$
\end{thm}

The sign vector $\tau$ acts by multiplying coordinate $\sigma_{ij}$ by
$\tau_{ij}$. 

\begin{proof}
The proof idea is to show that the toric ideals $I_{T,O}$ have the
same generators as the toric initial ideals of the Grassmannian that
have already been characterized in \cite{Speyer2004}.   
Without loss of generality, we may suppose that the leaves of $T$ are
labeled by $[n]$, that the tree is drawn without edge crossings,
and the leaves are labeled in increasing order from left to right.
These assumptions will allow us to ignore the sign vector $\tau$ in
the proof.  The sign vector results from straightening the tree 
and permuting the columns in the Steifel coordinates.  This 
results in sign changes in the Pl\"ucker coordinates.   

In Corollary \ref{thm:hidtree}, we saw that $I_{T,O}$ was generated by
the quadratic relations 
$$\sigma_{ik} \sigma_{jl} - \sigma_{il} \sigma_{kj}$$
such that there is a choke point in $T$ between $\{i,j\}$ and
$\{k,l\}$.  This is the same as saying that the induced subtree of $T$
on $\{i,j,k,l\}$ has the split $\{i,j\} |  \{k,l\}$.   These are
precisely the generators of the toric initial ideals of the
Grassmannian $G_{2,n}$ identified in \cite{Speyer2004}. 
\end{proof}

In the preceding Theorem, any weight vector $\omega$ that belongs to
the relative interior of the cone of the tropical Grassmannian
corresponding to the tree $T$ will serve as the desired partial term
order.  We refer to  \cite{Speyer2004} for background on the tropical
Grassmannian and toric degenerations of the Grassmannian.  Since and
ideal and its initial ideals have the same Hilbert function, we see
Catalan numbers emerging as degrees of Bayesian networks yet again.

\begin{cor}
Let $T$ be a rooted, directed, binary tree and $O$ consist of the
leaves of $T$.  Then $\deg I_{T,O} = \frac{1}{n-1} { 2n -4 \choose
  n-2}$, the $(n-2)$-nd Catalan number. 
\end{cor}

The fact that binary hidden tree models are toric degenerations of the
Grassmannian has potential use in phylogenetics.  Namely, it suggests
a family of new models, of the same dimension as the binary tree models,
that could be used to interpolate between the various tree models.
That is, rather than choosing a weight vector in a full dimensional
cone of the tropical Grassmannian, we could choose a weight vector
$\omega$ that sits inside of lower dimensional cone.  The varieties of
the initial ideals $V({\rm in}_\omega(I_{2,n}))$ then
correspond to models that sit somewhere ``between'' models
corresponding of the full dimensional trees of the maximal dimensional
cones containing $\omega$.  Phylogenetic recovery algorithms could
reference these in-between models to indicate some uncertainty about
the relationships between a given collection of species or on
a given branch of the tree.  These new models have the advantage that
they have the same dimension as the tree models and so there is no
need for dimension penalization in model selection.


\subsection{Matrix Schubert Varieties}

In this section, we will describe how certain varieties called matrix
Schubert varieties arise as special cases of the varieties of 
hidden variable models for Gaussian Bayesian networks.  More
precisely, the variety for the Gaussian Bayesian network will be the
cone over one of these matrix Schubert varieties.    To do this, we
first need to recall some equivalent definitions of matrix Schubert
varieties. 

Let $w$ be a partial permutation matrix, which is an $n \times n$
$0/1$ matrix with at most one $1$ in each row and column. 
The matrix $w$ is in the affine space $\cc^{n \times n}$.  The Borel
group $B$ of upper triangular matrices acts on $\cc^{n \times n}$ on
the right by multiplication and on the left by multiplication by the
transpose.    

\begin{defn}
The matrix Schubert variety $X_w$ is the orbit closure of $w$ by the
action of $B$ on the right and left: 
$$X_w =  \overline{B^TwB}.$$  
Let $I_w$ be the vanishing ideal of $X_w$.
\end{defn}

The matrix Schubert variety $X_w \subseteq \cc^{n\times n}$, so we can
identify its coordinate ring with a quotient of $\cc[\sigma_{ij} \, \,
  | \,\, i \in [n], j \in [n']]$.  Throughout this section
$[n'] = \{1', 2', \ldots, n'\}$, is a set of $n$ symbols that we use
to distinguish from $[n] = \{1,2, \ldots, n\}$.

An equivalent definition of a matrix
Schubert variety comes as follows.  Let $S(w) = \{ (i,j) \, | \, \,
w_{ij} =1 \}$ be the index set of the ones in $w$.  For each $(i,j)$
let $M_{ij}$ be the variety of rank one matrices: 
$$M_{ij} =  \left\{ x \in \cc^{n \times n} \, \, | \, \, {\rm rank}  \, x \leq
1,  x_{kl} = 0 \mbox{ if } k <i \mbox{ or } l < j \right \}.$$ 
Then 
$$X_w  =  \sum_{(i,j) \in S(w)}  M_{ij}$$
where the sum denotes the pointwise Minkowski sum of the varieties.
Since $M_{ij}$ are cones over projective varieties, this is the same
as taking the join, defined in the next section.

\begin{ex}
Let $w$ be the partial permutation matrix
$$w = \begin{pmatrix}
1 & 0 & 0 \\
0 & 1 & 0 \\
0 & 0 & 0 \end{pmatrix}.$$
Then $X_w$ consists of all $3 \times 3$ matrices of rank $\leq 2$ and
$I_w = \left< |\Sigma_{[3], [3']}| \right>$.
More generally, if $w$ is a partial permutation matrix of the form
$$w =  \begin{pmatrix}
E_d & 0 \\
0 & 0  \end{pmatrix},
$$
where $E_d$ is a $d \times d$ identity matrix, then $I_w$ is the ideal
of $(d+1)$ minors of a generic matrix.  \qed
\end{ex} 

The particular Bayesian networks which yield the desired varieties
come from taking certain partitions of the variables.  In particular,
we assume that the observed variables come in two types labeled by
$[n] = \{1,2, \ldots, n\}$ and $[n'] = \{1',2', \ldots, n' \}$.  The
hidden variables will be labeled by the set $S(w)$. 

Define the graph $G(w)$ with vertex set $V = [n] \cup [n'] \cup S(w)$
and edge set consisting of edges $k \to l$ for all $k < l \in [n]$,
$k' \to l'$ for all $k' < l' \in [n']$,  $(i,j) \to k$ for all $(i,j)
\in S(w)$ and $k \geq i$ and $(i,j) \to k'$ for all $(i,j) \in S(w)$
and $k' \geq j$. 

\begin{thm}\label{thm:schubert}
The generators of the ideal $I_w$ defining the matrix Schubert variety
$X_w$ are the same as the generators of the ideal $I_{G(w), [n] \cup
  [n']}$ of the hidden variable Bayesian network for the DAG $G(w)$
with observed variables $[n] \cup [n']$.  That is, 
$$I_w \cdot \cc[\sigma_{ij}\, \,  | \, \, i,j \in [n] \cup [n'] ]  =
I_{G(w), [n] \cup [n']}.$$ 
\end{thm}

\begin{proof}
The proof proceeds in a few steps.  First, we give a parametrization
of a cone over the matrix Schubert variety, whose ideal is naturally
seen to be $I_w \cdot \cc[\sigma_{ij}\, \,  | \, \, i,j \in [n] \cup
  [n'] ].$
Then we describe a rational transformation $\phi$ on $\cc[\sigma_{ij} \, \, |
  \, \, i,j \in [n] \cup [n'] ]$ such that $\phi (I_w) = I_{G(w), [n]
  \cup [n']}$.  We then exploit that fact that this transformation is
invertible and the elimination ideal $I_{G(w), [n] \cup [n']} \cap
\cc[\sigma_{ij} \,\, | \, \, i \in [n], j \in [n'] ]$ is fixed to deduce the desired equality.

First of all, we give our parametrization of the ideal $I_w$.  To do
this, we need to carefully identify all parameters involved in the
representation.  First of all, we split the indeterminates in the ring
$\cc[\sigma_{ij} \, \, | \, \, i,i \in [n] \cup [n']]$ into three classes of indeterminates:  those with
$i,j \in [n]$, those with $i,j \in [n']$, and those with $i \in [n]$
and $j \in [n']$.  Then we define a parametrization $\phi_{w}$ which is
determined as follows:

$$\phi_w: \cc[\tau, \gamma, a, \lambda] \rightarrow \cc[\sigma_{ij} \,
  \, | \, \, i,j \in [n] \cup [n']$$
$$\phi_w(\sigma_{ij})  = \left\{  
\begin{array}{ll}
\tau_{ij} &  \mbox{ if }  i,j \in [n] \\
\gamma_{ij} &  \mbox{ if } i,j \in [n'] \\
\sum_{ (k,l) \in S(w): k \leq i , l \leq j }  a_{(k,l)} \lambda_{(k,l),i}
\lambda_{(k,l),j} &  \mbox{ if }  i \in [n], j \in [n']
\end{array}  \right.$$

Let $J_w =  \ker \phi_w$.  Since the $\tau$, $\gamma$, $\lambda$, and
$a$ parameters are all algebraically independent, we deduce that in
$J_w$, there will be no generators that involve combinations of the
three types of indeterminates in $\cc[\sigma_{ij} \, \, | \, \, i,j
  \in [n] \cup [n']]$.  Furthermore,
restricting to the first two types of indeterminates, there will not
be any nontrivial relations involving these types of indeterminates.
Thus, to determine $J_w$, it suffices to restrict to the ideal among
the indeterminates of the form $\sigma_{ij}$ such that $i \in [n]$ and
$j \in [n']$.  However, considering the parametrization in this case,
we see that this is precisely the parametrization of the ideal $I_w$,
given as the Minkowski sum of rank one matrices.  Thus, $J_w = I_w$.

Now we will define a map from $\phi: \cc[\sigma_{ij}] \rightarrow
\cc[\sigma_{ij}]$ which sends $J_w$ to another ideal, closely related
to $I_{G(w), [n] \cup [n']} $.  To
define this map, first, we use the fact that from the submatrix
$\Sigma_{[n], [n]}$ we can recover the $\lambda_{ij}$ and $a_i$
parameters associated to $[n]$, when only considering the complete
subgraph associated to graph $G(w)_{[n]}$ (and ignoring the treks that
involve the vertices $(k,l) \in S(w)$).  This follows because these parameters
are identifiable by Proposition \ref{prop:dim}.  A similar fact holds
when restricting to the subgraph $G(w)_{[n']}$.  The ideal $J_w$ we
have defined thus far can be considered as the vanishing ideal of a
parametrization which gives the complete graph parametrization for
$G(w)_{[n]}$ and $G(w)_{[n']}$ and a parameterization of the matrix
Schubert variety $X_w$ on the $\sigma_{ij}$ with $i \in [n]$ and $j
\in [n']$.  So we can rationally recover the $\lambda$ and $a$
parameters associated to the subgraphs $G(w)_{[n]}$ and
$G(w)_{[n']}$.

For each $j<k$ pair in $[n]$ or in $[n']$, define the partial trek
polynomial 
$$s_{jk}(\lambda) =  \sum_{m =1}^{k-j} \sum_{j=l_0<l_1 <
    \ldots < l_m = k}  \prod_{i=1}^m  \lambda_{l_{i-1}l_{i}}. 
$$
We fit these into two upper triangular matrices $S$ and $S'$ where
$S_{jk} = s_{jk}$ if $j < k$ with $j,k \in [n]$, $S_{jj} = 1$ and
$S_{jk} = 0$ otherwise, with a similar definition for $S'$ with $[n]$
replaced by $[n']$.  Now we are ready to define our map.  Let $\phi$
be the rational map $\phi: \cc[\Sigma] \to \cc[\Sigma]$ which leaves
$\sigma_{ij}$ fixed if $i,j \in [n]$ or $i,j \in [n']$, and maps
$\sigma_{ij}$ with $i \in [n]$ and $j \in [n']$ by sending

$$\Sigma_{[n], [n']}  \mapsto   S  \Sigma_{[n],[n']} (S')^T.$$

This is actually a rational map, because the $\lambda_{ij}$ that
appear in the formula for $s_{jk}$ are expressed as rational functions
in terms of the $\sigma_{ij}$ by the rational parameter recovery
formula of Proposition \ref{prop:dim}.  Since this map transforms
$\Sigma_{[n], [n']}$ by multiplying on the left and right but lower
and upper triangular matrices, this leaves the ideal $J_w \cap
\cc[\sigma_{ij} \, \, | \,\, i \in [n], j \in [n']]$ fixed.  Thus $J_w
\subseteq \phi(J_w)$.  On the other hand $\phi$ is invertible on $J_w$
so $J_w = \phi(J_w)$.

If we think about the formulas for the image $\phi \circ
\phi_w$, we see that the formulas for $\sigma_{ij}$ with $i \in [n]$
and $j \in [n']$ in terms of parameters are the correct
formulas which we would see coming from the parametrization
$\phi_{G(w)}$.  On the other hand, the formulas for $\sigma_{ij}$ with
$i,j \in [n]$ or $i,j \in [n']$ are the formulas for the restricted
graph $G_{[n]}$ and $G_{[n']}$, respectively.  Since every trek
contained in $G_{[n]}$ or $G_{[n']}$ is a trek in $G(w)$, we see that
the current parametrization of $J_w$ is only ``almost correct'', in
that it is only missing terms corresponding to treks that go outside
of $G(w)_{[n]}$ or $G(w)_{[n']}$.  Denote this map by $\psi_w$, and
let $\phi_{G(w)}$ be the actual parametrizing map of the model.  Thus,
we have, for each $\sigma_{ij}$ with $i,j \in [n]$ or $i,j \in [n']$,
$\phi_{G(w)}(\sigma_{ij}) =  \psi_w(\sigma_{ij})  +
r_w(\sigma_{ij})$, where $r_w(\sigma_{ij})$ is a polynomial remainder
term that
does not contain any $a_i$ with $i \in [n] \cup [n']$, when $i,j \in
[n]$ or $i,j \in [n']$, and $r_w(\sigma_{ij}) = 0$ otherwise.  On the other
hand, every term of $\psi_w(\sigma_{ij})$ will involve exactly one
$a_i$ with $i \in [n] \cup [n']$, when $i,j \in [n]$ or $i,j \in
[n']$.  

Now we define a weight ordering $\prec$ on the ring $\cc[a,\lambda]$ that
gives $\deg a_i = 1$ if $ i \in [n] \cup [n']$ and $\deg a_i = 0$
otherwise and $\deg \lambda_{ij} = 0$ for all $i,j$.  Then, the
largest degree term of $\phi_{G(w)}(\sigma_{ij})$ with respect to this
weight ordering is $\psi_w(\sigma)$.  Since highest weight terms must all
cancel with each other, we see that
$f \in I_{G(w), [n] \cup [n']}$, implies that  $f \in J_w$.  Thus, we
deduce that $I_{G(w), [n] \cup [n']}  \subseteq J_w$.  On the other
hand, 
$$I_{G(w), [n] \cup [n]'}  \cap \cc[\sigma_{ij} \, \, | \, \, i \in
  [n], j \in [n']  ]  = J_w \cap \cc[\sigma_{ij} \, \,| \, \, i \in
  [n], j \in [n']  ]$$
and since the generators of $J_w \cap \cc[\sigma_{ij} \, \,| \, \, i \in
  [n], j \in [n']  ]$ generate $J_w$, we deduce that $J_w \subseteq
  I_{G(w), [n] \cup [n']}$ which completes the proof.
\end{proof}

The significance of Theorem \ref{thm:schubert} comes from the work of
Knutson and Miller \cite{Knutson2005}.  They gave a complete
description of antidiagonal Gr\"obner bases for the ideals $I_w$.
Indeed, these ideals are generated by certain subdeterminants of the 
matrix $\Sigma_{[n],[n']}.$  These determinants can be interpretted
combinatorially in terms of the graph $G(w)$.

\begin{thm}\cite{Knutson2005}
The ideal $I_w$ defining the matrix Schubert variety is generated by
the conditional independence statements implied by the DAG $G(w)$.  In
particular, 
$$I_w  =  \left<  \#C+1 \mbox{ minors } of \Sigma_{A,B} \, \,  | \,\, A
\subset [n], B \subset [n'], C \subset S(w), \mbox{ and }  C \mbox{
  $d$-separates } A \mbox{ from } B \right>.$$
\end{thm}


\subsection{Joins and Secant Varieties}

In this section, we will show how joins and secant varieties arise as
special cases of Gaussian Bayesian networks in the hidden variable
case.  This, in turn, implies that techniques that have been developed
for studying defining equations of joins and secant varieties (e.g.
\cite{Landsberg2004, Sturmfels2006}) might be useful for studying the
equations defining these hidden variable models. 

Given two ideals $I$ and $J$ in a polynomial ring  $\kk[\bfx] =
\kk[x_1, \ldots, x_m]$, their \emph{join} is the new ideal 
$$I*J  :=    \left(  I(\bfy)  + J(\bfz)  +  \left<  x_i - y_i - z_i \,
\, | \, \, i \in [m] \right>  \right)  \bigcap  \cc[\bfx]$$ 
where $I(\bfy)$ is the ideal obtained from $I$ by plugging in the
variables $y_1, \ldots, y_m$ for $x_1, \ldots, x_m$.  The secant ideal
is the iterated join: 
$$I^{\{r\}}  =  I \ast I \ast \cdots \ast I$$
with $r$ copies of $I$.  If $I$ and $J$ are homogeneous radical ideals
over an algebraically closed field, the join ideal $I \ast J$ is the
vanishing ideal of the join variety which is defined geometrically by
the rule 
$$V(I \ast J)  =  V(I) \ast V(J)  =    \overline{\bigcup_{a \in V(I)}
\bigcup_{b \in V(J)}   < a,b>}$$ 
where $<a,b>$ denotes the line spanned by $a$ and $b$ and the bar denotes the Zariski closure. 

Suppose further that $I$ and $J$ are the vanishing ideals of
parametrizations;  that is there are $\phi$ and $\psi$ such that 
$$\phi:  \cc[\bfx] \rightarrow  \cc[\theta]  \mbox{ and }  \psi:
\cc[\bfx]  \rightarrow  \cc[\eta]$$   
and $I  = \ker \phi$ and $J = \ker \psi$.  Then  $I \ast J$ is the
kernel of the map 
$$\phi + \psi  :  \cc[\bfx]  \rightarrow  \cc[\theta, \eta]$$
$$    x_i  \mapsto   \phi(x_i)  +  \psi(x_i).$$

Given a DAG $G$ and a subset $K \subset V(G)$,  $G_K$ denotes the
induced subgraph on $K$. 

\begin{prop}
Let $G$ be a DAG and suppose that the vertices of $G$ are partitioned
into $V(G)  = O \cup H_1 \cup H_2$ where both $H_1$ and $H_2$ are
hidden sets of variables.  Suppose further that there are no edges of
the form $o_1 \to o_2$ such that $o_1, o_2 \in O$ or edges of the form
$h_1 \to h_2$ or $h_2 \to h_1$ with $h_1 \in H_1$ and $h_2 \in H_2$.
Then 
$$I_{G,O}  =    I_{G_{O \cup H_1}, O}  \ast  I_{G_{O \cup H_2}, O}.$$ 
\end{prop}

The proposition says that if the hidden variables are partitioned with
no edges between the two sets and there are no edges between the
observed variables the ideal is a join. 

\begin{proof}
The parametrization of the hidden variable model only involves the
$\sigma_{ij}$ such that $i,j \in O$.  First, we restrict to the case
where $i \neq j$.  Since there are no edges
between observed variables and no edges between $H_1$ and $H_2$, every
trek from $i$ to $j$ involves only edges in $G_{O \cup H_1}$ or only
edges in $G_{O \cup H_2}$. 
This means that 
$$\phi_G(\sigma_{ij})  =    \phi_{G_{O \cup H_1}}(\sigma_{ij})  +
\phi_{G_{O \cup H_2}}(\sigma_{ij})$$ 
and these summands are in non-overlapping sets of indeterminates.
Thus, by the comments preceding the proposition, the ideal only in the
$\sigma_{ij}$ with $i \neq j \in O$ is clearly a join.  However, the
structure of this hidden variable model implies that there are no
nontrivial relations that involve the diagonal elements $\sigma_{ii}$
with $i \in O$.  This implies that $I_{G,O}$ is a join. 
\end{proof}

\begin{ex}
Let $K_{p,m}$ be the directed complete bipartite graph with
bipartition $H = [p']$ and $O = [m]$ such that $i' \to j \in
E(K_{p,m})$ for all $i' \in [p']$ and $j \in [m]$.  Then
$K_{p,m}$ satisfies the conditions of the theorem recursively up to
$p$ copies, and we see that: 
$$I_{K_{p,m}, O}   =   I_{K_{1,m}, O}^{\{p\}}.$$
This particular hidden variable Gaussian Bayesian network is known as
the factor analysis model.  This realization of the factor analysis model as a
secant variety was studied extensively in \cite{Drton2006}. 

\end{ex}

\begin{ex}
Consider the two ``doubled trees'' pictured in the figure.  
\begin{center}
\resizebox{4cm}{!}{\includegraphics{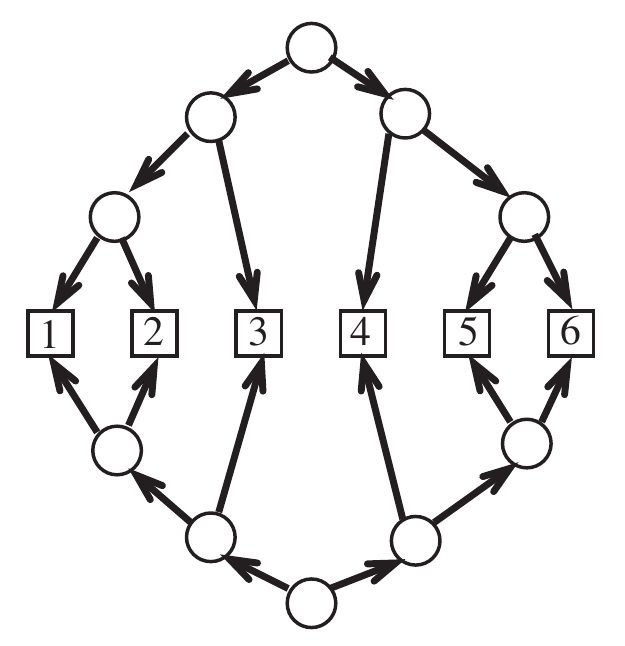}} \hspace{1cm}
\resizebox{4cm}{!}{\includegraphics{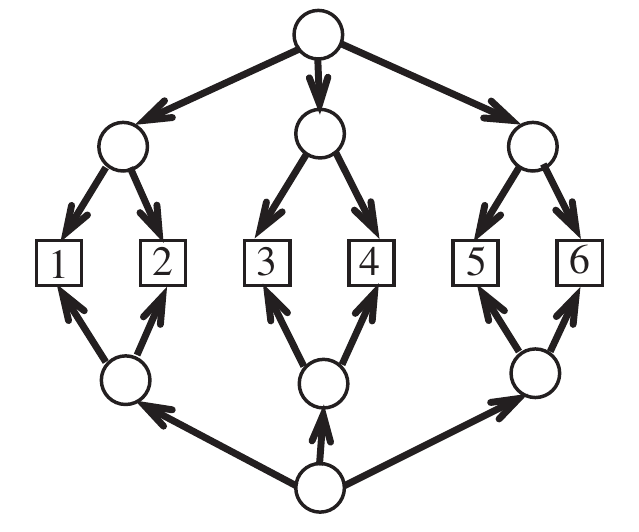}}
\end{center}

Since in each case, the two subgraphs $G_{O \cup H_1}$ and $G_{O \cup
H_2}$ are isomorphic, the ideals are secant ideals of
the hidden tree models $I_{T,O}$ for the appropriate underlying trees.  In
both cases, the ideal $I_{T,O}^{\{2\}}  =  I_{G,O}$ is a principal ideal,
generated by a single cubic.  In the first case, the ideal is the
determinantal ideal $J_T^{\{2\}}  =  \left<  |  \Sigma_{123,456} |
\right> $.  In the second case, the ideal is generated by an eight
term cubic 

\begin{eqnarray*}
\hspace{3cm}  I_{G,O} & = &  \left<  \sigma_{13} \sigma_{25}
 \sigma_{46}  -  \sigma_{13}  \sigma_{26} \sigma_{45}  -
 \sigma_{14}\sigma_{25} \sigma_{36} +  \sigma_{14}  \sigma_{26}
 \sigma_{35}  \right. \\ 
 &  &   \left.   +  \sigma_{15}  \sigma_{23}  \sigma_{46}
  -  \sigma_{15} \sigma_{24}  \sigma_{36}   -  \sigma_{16} \sigma_{23}
 \sigma_{45} +  \sigma_{16}  \sigma_{24} \sigma_{35} \right>. 
\end{eqnarray*}
\qed
\end{ex}

In both of the cubic cases in the previous example, the ideals under
questions were secant ideals of toric ideals that were initial ideals
of the Grassmann-Pl\"ucker ideal, as we saw in Theorem \ref{thm:grin}.
Note also that the secant ideals $I_{T,O}^{\{2\}}$ are, in fact, the
initial terms of the $6 \times 6$ Pfaffian with respect to
appropriate weight vectors.  We conjecture that this pattern
holds in general.

\begin{conj}
Let $T$ be a binary tree with $n$ leaves and $O$ the set of leaves of
$T$.   Let $I_{2,n}$ be the
Grassmann-Plu\"ucker ideal, let $\omega$ be a weight vector and $\tau$
a sign vector so that $I_{T,O}  =  \tau \cdot {\rm in}_\omega(I_{2,n})$ as
in Theorem \ref{thm:grin}.  Then for each $r$ 
$$I_{T,O}^{\{r\}}  =   \tau \cdot {\rm  in}_\omega(I^{\{r\}}_{2,n}).$$ 
\end{conj}


\end{document}